\documentclass[11pt]{article}                                                 
\usepackage{amssymb,chicago, epsf, here, amsfonts,latexsym,dsfont,amsmath}       

\usepackage{color}
\usepackage{ulem}

\advance \topmargin by -\headheight                                           
\advance \topmargin by -\headsep                                              
\advance \topmargin .2in                                                      
\newcommand{\RR}{\mathbb {R}}
\newcommand{\PP}{\mathbb {P}}
\newcommand{\EE}{\mathbb {E}}

\def\ben{\begin{eqnarray}}
\def\be*{\begin{eqnarray*}}
\def\non{\end{eqnarray}}
\def\no*{\end{eqnarray*}}
\def\argmax{\mathop{\rm argmax}}

\newtheorem{theorem}{Theorem}[section]
\newtheorem{lemma}[theorem]{Lemma}

\begin{document}
\title{Removing the Curse of Superefficiency: an Effective Strategy For Distributed Computing in Isotonic Regression}
\author{Moulinath Banerjee$^1$ and Cecile Durot$^2$}
\date{} 
\maketitle
{\small\begin{center}
$^{1}$University of Michigan,
451, West Hall, 
1085 South University 
Ann Arbor, MI 48109\\
$^{2}$Modal'x, Universit\'e Paris Nanterre, F-92001, Nanterre, France\\
\end{center}}
\begin{abstract}
We propose a strategy for computing the isotonic least-squares estimate of a monotone function in a general regression setting where the data are distributed across different servers and the observations across servers, though independent, can come from heterogeneous sub-populations, thereby violating the identically distributed assumption. Our strategy fixes the super-efficiency phenomenon observed in prior work on distributed computing in %non-standard problems featuring non-Gaussian limits and slow convergence rates, and in particular
the isotonic regression framework, where averaging several isotonic estimates (each computed at a local server) on a central server produces super-efficient estimates that do not replicate the properties of the global isotonic estimator, i.e. the isotonic estimate that would be constructed by transferring \emph{all the data} to a single server. The new estimator proposed in this paper works by smoothing the data on each local server, communicating the smoothed summaries to the central server, and then computing an isotonic estimate at the central server, and is shown to replicate the asymptotic properties of the global estimator, and also overcome the super-efficiency phenomenon exhibited by earlier estimators. For data on $N$ observations, the new estimator can be constructed by transferring data just over order $N^{1/3}$ across servers [as compared to transferring data of order $N$ to compute the global isotonic estimator], and requires the same order of computing time as the global estimator. 


\end{abstract} 

\section{Background} 
Distributed computing has now become significant in the practice of statistics as well as other branches of data science. Large volumes of data, often relating 
to the same or closely related studies or experiments, are no longer stored on one single computer; rather, they are distributed across a number of platforms in 
some structured manner, owing partly to natural memory constraints on individual machines, and partly for convenience. This, typically, poses problems for computing optimal estimates of parameters of interest from the data at hand. Conventional statistical estimates are generally obtained under the premise that the totality of the data is accessible to a single computing device and can be processed at one stroke, yielding estimates that are \emph{optimal} in some quantitatively defined sense. However, this is not automatically the case in a distributed environment. The calculation of \emph{global estimates} that require simultaneous processing of all available data then entails transferring the entire bulk of data from different computers to a central machine, which in itself can 
be both time and resource consuming, followed by a potentially complex computation on the aggregated data (of massive volume), which may be infeasible under many circumstances. 
\newline
\newline
Divide and conquer algorithms are a standard approach to addressing these issues in a distributed computing environment. The idea behind this is as follows: suppose the entire data set is stored across a number of machines. On each machine, calculate a natural estimate of the parameter of interest from the data on it and transfer this estimate to a central machine. Next, combine the estimators thus obtained, at the central machine in a judicious way to produce a final 
estimate, the so-called \emph{pooled estimate}, which replicates the properties of the natural \emph{global estimate}, i.e. the one we could have computed were it feasible to store and analyze all available data on one machine. The term `replicates the properties' can be understood in various ways and is often specific to the problem at hand: one might be able to show that the pooled and the global estimates have the same rate of convergence, or are comparable, up to constants, in terms of a certain measure of risk, or it might even be possible to demonstrate that the pooled estimate and the global estimate have the same limit distributions under appropriate conditions. The other important factor is computational burden: one would expect that the divide and conquer algorithm is not substantially more computationally onerous than the global estimator.  As the literature on distributed computing is enormous, here we provide a selection of  instances of research on distributed computing problems in a variety of statistical/machine-learning contexts: see, e.g. ~\cite{hsieh2014divide}, ~\cite{LiEtAl13}, ~\cite{Zhang13}, ~\cite{Zhao14}, ~\cite{battey2015distributed}, ~\cite{shang2017computational}, ~\cite{volgushev2017distributed}. The above papers illustrate that the sample splitting approach buys computational dividends, yet statistical optimality [in the sense that the resulting estimator is as efficient (or minimax rate optimal) as the global estimate based on applying the estimation algorithm to the entire data set] is retained. 
\newline
\newline
In the nonparametric function estimation context, most results of the divide and conquer type focus on estimation under smoothness constraints, where the essence of the strategy is to compute a smoothed estimator of the unknown function at each server and combine the estimators at the central server, by averaging; this strategy is employed, for example, in \cite{LiEtAl13}, ~\cite{Zhang13}, ~\cite{Zhao14}. However, the averaging strategy leads to highly problematic pooled estimators in non-regular function estimation problems, e.g. function estimation under a monotonicity constraint, where the least squares estimates under the monotonicity constraint are non-standard/non-regular in the sense that they are highly non-linear in the data, and exhibit non-Gaussian limit distributions. This is the core content of the recent work by \cite{banerjee2018divide} [henceforth BDS] where it is demonstrated that in monotone function estimation problems, the `pooled-by-averaging' estimator [henceforth, generally referred to as BDSE] becomes \emph{super-efficient}: its asymptotic relative efficiency (in terms of MSE) with respect to the global monotone least squares estimator computed at any single model goes to infinity, whereas, in the uniform sense, the ARE goes to 0, i.e. the maximal MSE of BDSE over a collection of models relative to that of the global least squares estimator goes to $\infty$. Furthermore, BDSE has a distribution different from that of the global estimator which converges to a Chernoff limit (discussed in details below). Indeed, the super-efficiency property just alluded to, has also been observed for the pooled-by-averaging estimator in the genre of non-standard problems exhibiting cube-root asymptotics in the sense of ~\cite{KP90}; see, ~\cite{shi2017massive}. 
\newline
\newline
Our goal in this paper is to construct a \emph{new estimator} in the monotone function estimation problem which does not suffer from the super-efficiency problem and which also exhibits the limiting properties of the global estimator. To this end, we provide some details of the problem considered in BDS and the results obtained, as they are crucial to understanding the goal and the approach of the current work.  
\newline
\newline
Consider a sample of size $N$ (very large) from the model $Y_i = \mu(X_i) + \epsilon_i$ which is distributed across $m$ different servers, each server containing a subsample of size $n$, and $m = o(N)$. The function $\mu$ is known to be monotone and the $X_i$ come from a density on $[0,1]$. The residual 
$\epsilon_i$ is assumed to satisfy $E(\epsilon_i|X_i) = 0$. Computing the global isotonic estimate at a point $t_0 \in (0,1)$ involves moving all the data to a central server and performing the isotonization on all $N$ data-points on the central server. This can be time-consuming when $N$ is really large. The construction of BDSE involves computing the isotonic estimate of $\mu$, say $\widehat{\mu}_j$, on the $j$'th server, and then obtaining the average of these isotonic estimates. Hence, the pooled estimate at the point $t_0$ is given by:  $\overline{\mu}(t_0) := m^{-1}\,\sum_{j=1}^m\,\widehat{\mu}_j(t_0)$. Computing BDSE at a particular point only requires transferring $m$ numbers (from the $m$ machines) to the central server, where $m = o(N)$. 
\newline
\newline
One can compare the computational burden involved in the calculation of the global estimator to that for BDSE. For the global estimator, once all the data-points have been transferred to the central machine, sorting of the $X_i$'s (resulting in an induced sorting of the $Y_i$'s) can be accomplished typically in $O(N \log N)$ time. Post-sorting, one can implement isotonic regression via the PAVA algorithm \cite{RWD88} (Chapter 1) which takes $O(N)$ time. Thus, the total computational burden is $O(N \log N)$ computing time plus the transferring of $N$ bivariate pairs to the central machine. On the other hand, for the pooled estimator, on each machine, the isotonic estimate based on the subsample stored in that machine takes $O(n \log n)$ computing time,  leading to a total computing time of $O( m n \log n)$. At the central server, averaging takes $O(m)$ time. If $n \sim N^{\gamma}$ for some $0 < \gamma < 1$, this gives a total computing time of order $O(N \log N)$, and in addition, one transfers $m \sim N^{1-\gamma}$ scalars (the values $\widehat{\mu}_j(t_0)$ for $j = 1, 2, \ldots, m$) to the central machine. Thus, the pooled estimator is computationally less burdensome than the global estimator. Similar considerations apply to the computation of the global and pooled isotonic estimators of the inverse function $\mu^{-1}$. 
\newline
BDS showed that their pooled-by-averaging estimator (BDSE) of the inverse function has \emph{dichotomous behavior}. We briefly revisit this important result. For convenience and the sake of completeness, we state these results essentially in their entirety. 
\newline 
Consider a {nonincreasing} and continuously differentiable function $\mu_0$ on $[0,1]$ with $0 < c < |\mu_0^{'}(t)| < d < \infty$ for all $t \in [0,1]$.  For an $x_0 \in (0,1)$, define a neighborhood $\mathcal{M}_0$ of $\mu_0$ as the class of all continuous nonincreasing functions $\mu$ on $[0,1]$  that are continuously differentiable on $[0,1]$, coincide with $\mu_0$ outside of $(x_0 - \epsilon_0, x_0 + \epsilon_0)$ for some (small) $\epsilon_0>0$,  satisfy $0 < c < |\mu'(t)| < d < \infty$ for all $t \in [0,1]$, and such that $\mu^{-1}(a) \in (x_0 - \epsilon_0, x_0 + \epsilon_0)$ where $a = \mu_0(x_0)$. Now, consider $N$ i.i.d.~observations $\{(X_i, Y_i)\}_{i=1}^N$ from {$(X,Y)$ where $Y_i = \mu_0(X_i) + \epsilon_i$
and $X_i \sim \mbox{Uniform}(0,1)$ is independent of $\epsilon_i \sim N(0, v^2)$. Then, the isotonic estimate $\hat{\theta}_N$ of $\theta_0: = 
\mu_0^{-1}(a)$ (which is $x_0$) satisfies
\begin{equation}\label{eq:IsoRegLimDist}
N^{1/3}\,(\hat{\theta}_N - \theta_0) \stackrel{d}{\rightarrow} G,
\end{equation}
as $N\to\infty$, where $G =_d \tilde{\kappa}\,\mathbb{Z}$, with $\mathbb{Z}$ following the Chernoff distribution, and $\tilde{\kappa} > 0$ being a constant. Writing $N= m \times n$, where $m$ and $n$ are as defined above, as $N \rightarrow \infty$, the BDSE of $\mu^{-1}(a)$, say $\overline{\theta}_m$ satisfies: 
$$N^{1/3}(\overline{\theta}_{m} - \theta_0) \rightarrow_d m^{-1/6} H\,,$$
where $H$ has the same variance as $G$ but is distributed differently from $G$. Furthermore, 
\begin{equation}\label{eq:MSEGlobalEst}
\mathbb{E}_{\mu_0}\left[N^{2/3}(\hat{\theta}_N - \theta_0)^2\right] \rightarrow \mathrm{Var}(G) \quad {\mbox{ and }}\quad
%\qquad \mbox{as } N \to \infty,
%\end{equation}
%while
%\begin{equation}\label{eq:MSEPooledEst} 
\mathbb{E}_{\mu_0}\left[N^{2/3}(\overline{\theta}_{m} - \theta_0)^2\right] \rightarrow m^{-1/3}\,\mathrm{Var}(G),
%\qquad \mbox{as } N \to \infty\,.
\end{equation}
as $N\to\infty$. Hence, %for estimating $\theta_0 = \mu_0^{-1}(a)$, 
BDSE \emph{outperforms} the global inverse isotonic regression estimator in terms of point wise MSE. 
\newline
\newline
This phenomenon is \emph{reversed} when one looks at the maximal MSEs of the two estimators over the class of models defined by $\mathcal{M}_0$, as 
described in Theorem 5.1 of BDS.
\begin{theorem}(Theorem 5.1 of BDS)\label{superefficient}
Let
\begin{equation}\label{eq:IsoRegUnifBd}
E := \limsup_{N \to \infty}\,\sup_{\mu \in \mathcal{M}_0}\mathbb{E}_{\mu}\left[N^{2/3}(\hat{\theta}_N - \mu^{-1}(a))^2 \right]
\;  {\mbox{ and }}\;   
E_m: =  \liminf_{N \to \infty}\,\sup_{\mu \in \mathcal{M}_0}\mathbb{E}_{\mu}\left[ N^{2/3}(\overline{\theta}_{m} - \mu^{-1}(a))^2 \right]\end{equation}
where the subscript $m$ indicates that the maximal risk of the $m$-fold pooled estimator ($m$ fixed) is being considered. Then $E < \infty$ while $E_m \geq m^{2/3}\,c_0,$ for some $c_0 > 0$. When $m = m_n$ diverges to infinity, 
$$\liminf_{N \to \infty}\,\sup_{\mu \in \mathcal{M}_0}\mathbb{E}_{\mu}\left[ N^{2/3}(\overline{\theta}_{m_n} - \mu^{-1}(a))^2 \right] = \infty \,.$$ 
\end{theorem} 
Therefore, from Theorem~\ref{superefficient} it follows that the asymptotic maximal risk of BDSE diverges to $\infty$ at rate (at least) $m^{2/3}$. Thus, the better off we are with BDSE for a fixed function, the worse off we are in the uniform sense over the class of functions~$\mathcal{M}_0$. Hence, unfortunately, while maintaining a computational burden that is better than the global estimator, BDSE has undesirable statistical properties as seen above. 
\newline
\newline
As mentioned above, we will construct a corrected estimator in the isotonic regression problem that does not suffer from the undesirable `super-efficiency phenomenon'. Our new estimator recovers all the desirable properties of the global isotonic estimator, and is computationally not anymore onerous than the latter. Furthermore, for our analysis, we address a \emph{much broader scenario than is conventionally considered in isotonic regression problems.} Since we are thinking of large $N$ problems, with the data being stored separately across different servers, it is natural to allow heterogeneity in data. Thus, while considering the pairs $\{X_i, Y_i\}_{i=1}^N$ to be independent, we will \emph{no longer consider them to be identically distributed}; rather, they will be assumed to come from a number of different ($m$) sub-populations with the $(X_i, Y_i)$ pairs in each sub-population being i.i.d. What links the different sub-populations is the common mean function $\mu$ of interest: $E(Y_i | X_i) = \mu(X_i)$ for all $i$, for a common monotone function $\mu$. Furthermore, the $N$ pairs will be scrambled across a number of different servers (say $L$), with the same server hosting data from different sub-populations, as well as data from the same sub-population potentially stored on multiple servers\footnote{In BDS, the number of servers was designated by $m$, while in this paper we change notation and call it $L$. As we will see below, it is the number of different sub-populations that really enters into the properties of the pooled estimator in general and \emph{not} the number of servers. When each sub-population has its own server, then obviously $L=m$.} . Our new estimator essentially reverses the steps involved in constructing BDSE. BDSE involves isotonization on local servers followed by smoothing (as in simple averaging) on the central server, while, in this paper, we do the opposite: first smooth (by local averaging) on each local server, and then isotonize the smoothed data on the central server.  

\section{The Set-Up and the Estimator}
Assume that we have $m$ samples of respective sizes $n_1,\dots,n_m$ and that for all $j=1,\dots,m$, the $j$-th sample is composed of i.i.d. pairs  of real valued random variables $(X_{ji},Y_{ji}),$ $i=1,\dots,n_j$, such that $E(Y_{ji}|X_{ji})=\mu(X_{ji})$ for all $i,j$ and an unknown regression function $\mu$ defined on $[0,1]$. We denote by $F_{Xj}$ the common distribution function of the covariates $X_{ji},\ i=1,\dots,n_j$ in the $j$-th sample. The data are stored on several servers numbered $1,\dots,L$ for some integer $L\geq 1$. The allocation of data on the different servers is arbitrary in the sense that a sample can be spread on several servers, a server can host data from several different samples, and the number of stored observations can vary across the different servers. The number $L$ of different servers can even grow as $N\to\infty$. The total sample size is 
$$N=\sum_{j=1}^mn_j.$$

\bigskip

For ease of exposition, when considering simultaneously all the samples, we relabel the observations from the $m$ samples to obtain independent pairs $(X_i,Y_i)$, $i=1,\dots,N$ such that $E(Y_i|X_i)=\mu(X_i)$, where the distribution function of $X_i$ is one of $F_{X1},\dots,F_{Xm}$.
Let $K$ be a positive integer that grows to infinity as $N\to\infty$, and for all $k\in\{1,\dots,K\}$, let $I_k=((k-1)/K,k/K]$. Let $S_{\ell}$ denote the set of indices $i$, such that $(X_i, Y_i)$ is stored in the $\ell$'th server. Now, for each server $\ell$ ($1 \leq \ell \leq L$) record 
$$T_{\ell k}=\sum_{i=1}^{N} Y_i \mathds{1}_{i\in S_\ell}\mathds{1}_{X_i\in I_k}
$$
and
$$C_{\ell k}=\sum_{i=1}^{N} \mathds{1}_{i\in S_\ell }\mathds{1}_{X_i\in I_k},
$$
for $k\in\{1,\dots,K\}$. Next, for each $\ell$, transfer $\{(T_{\ell k}, C_{\ell k})\}_{k=1}^K$ to a central server. Compute a regressogram estimate on the central server in the following manner:
  for each $k\in\{1,\dots,K\}$, 
  %\comC{in the definition of $\overline y_k$ below, there was a typo : the sums are over $l=1,\dots,L$}
\begin{eqnarray*}
\overline y_k&=&\frac1{\sum_{\ell=1}^ L C_{\ell k}}\sum_{\ell =1}^L T_{\ell k}\\
&=&\frac 1{\sum_{i=1}^N \mathds{1}_{X_i\in I_k}}\sum_{i=1}^N Y_i\mathds{1}_{X_i\in I_k}.
\end{eqnarray*}
Our final estimator of   $(\mu(\overline x_1),\dots,\mu(\overline x_K))^T$, where $\overline x_k=k/K$, is
\begin{equation}
\label{argmin-LS}
\widehat y=\arg\min_{h\in\RR^K: h_1\geq\dots\geq h_K}\sum_{k=1}^Kw_k(\overline y_k-h_k)^2\,,
\end{equation} 
where
$$w_k=\frac 1{N}\sum_{i=1}^N\mathds{1}_{X_i\in I_k} = \frac{\sum_{\ell=1}^L C_{\ell k}}{N} \,.$$
\emph{Note that the estimator does not depend on the way the observations were stored across different servers.}
\newline
%{\bf I've rewritten part of the computational considerations. Most of the new material is at the end. Please give it a read.} 
\newline
{\bf Computational considerations:} Consider the computational burden for the new estimator. Assume, for now, that $K \sim N^{\zeta}$ for some $0 \leq \zeta < 1$. First, focus on the computational time it takes for calculating $(T_{\ell k}, C_{\ell k})$ for all $\ell$ and $1 \leq k \leq K$. 
%There are $n_j$ pairs $(X_i, Y_i)$ on server $S_j$. 
For each $X_i$, one has to determine in which interval $I_k$ it falls, and then assign the pair $(X_i, Y_i)$ to the interval $I_k$. This can be accomplished in $O(\log N^{\zeta}) = O(\log N)$ time. Since there are $N$ such points (scrambled across the different servers), the total time taken is $O(N \log N)$. Next, computing $(C_{\ell k}, T_{\ell k})$ for a fixed $\ell$ involves less than $2\,n_{\ell k}$ additions, where $n_{\ell k}$ is the number of $(X_i, Y_i)$ pairs assigned to $I_k$ on server $\ell$. Hence, computing the vector $\{C_{\ell k}, T_{\ell k}\}_{1 \leq k \leq K}$ takes $O(\sum_k n_{\ell k})$ time.  
Summing up across the different $\ell$'s, we are looking at a total of $O(N \log N) \vee O(N)$ time, i.e. $O(N \log N)$ time. 
\newline
After the pairs $\{T_{\ell k}, C_{\ell k}\}_{1\leq k \leq K}$ have been transferred to the central server, computing the vector $\{(w_k,\overline{y}_k)\}_{1\leq k \leq K}$ takes $O(L N^{\zeta})$ time, and the final isotonization step takes $O(N^{\zeta})$ time. Thus, the total computing time is $O(L N^{\zeta}) \vee 
O(N \log N)$ which is dominated by $O(N \log N)$ provided $L$ (which could grow with $N$) and $\zeta$ are not too large. In addition to the total computing time, the burden also involves transferring about $2LK \sim L N^{\zeta}$ numbers between machines, which is larger than the amount of data transferred in the construction of BDSE. As shall be seen below, with $K$ slightly larger than $N^{1/3}$ -- say $K \sim N^{1/3\,+\eta_1}$ ($\eta_1$ small) -- and $m$ of a smaller order than $N^{1/3}$, the new estimator is able to recover the properties of the global estimator: hence, so long as the number of machines is not too large --  say $L = N^{1/3\,-\eta_2}$ -- the total amount of data required to be transferred is of order $N^{2/3 \,+ \eta_1 - \eta_2} = o(N^{2/3})$ when $\eta_2 > \eta_1$. 
\newline
Note that the computation of the global isotonic estimator in this situation would require transferring all data points to the central server which is exactly 
$O(N)$ and the isotonic algorithm at the central server would take $O(N \log N)$ time. Note also that the minimum amount of data transferring needed for the 
new estimator above is of order $K$ (this happens when the number of servers $L$ is held fixed) and therefore of \emph{larger order than $N^{1/3}$}. On the other hand, in the scenario of BDS, where $L=m$, the BDSE is constructed using $m$ sub-samples where $m$ 
is of order at most $N^{1/4}$: this corresponds to a data-transfer of \emph{order at most $N^{1/4}$} numbers to construct the super-efficient estimator at any given point. The additional amount of data that needs to be transferred to construct the new estimator can be viewed as the cost of alleviating the super-efficiency phenomenon exhibited by BDSE. 
\newline
\newline
{\bf Characterization of the new estimator:} It is a standard result in isotonic regression that the minimum in \eqref{argmin-LS} is achieved at a unique vector $(\widehat y_{1},\dots,\widehat y_{K})^T$. We give below a characterization of the minimizer. 
%Let  $\overline x_{(1)}<\dots<\overline x_{(K)}$ denote the order statistics corresponding to $\overline x_1,\dots,\overline x_K$, and $\overline y_{(i)}$ the averaged observation corresponding to $\overline x_{(i)}$. 
In the sequel, we consider the piecewise-constant left-continuous estimator $\widehat\mu_N$ that is constant on  the intervals $[0, \overline x_{1}]$, and $(\overline x_{k-1},\overline x_{k}]$ for all $k=2,\dots,K$, and such that 
$$\widehat\mu_N(\overline x_{k})=\widehat y_{k}$$ 
for all $k=1,\dots,K.$ 
Let $F_N$ be the empirical distribution function corresponding to $X_1,\dots,X_N$ 
\begin{equation}
\label{def: FN}
F_N(x)=\frac 1N\sum_{i=1}^N\mathds{1}_{X_i\leq x},\ x\in\RR,
\end{equation}
and let 
$\Lambda_N$ be the piecewise-constant right-continuous process on $[0,1]$ that is constant on  the intervals $[0, \overline x_{1})$, and $[\overline x_{k-1},\overline x_{k})$ or all $k=2,\dots,K$ such that 
$$\Lambda_N\left(\overline x_j\right)=\sum_{k=1}^jw_k\overline y_{k}=\frac 1N\sum_{i=1}^NY_i\mathds{1}_{X_i\leq\overline x_j}$$
for all $j=1,\dots,K$, and $\Lambda_N(0)=0$.
Then, 
$$F_N\left(\overline x_j\right)=\sum_{k=1}^jw_k$$
and $\widehat \mu_N$  is the left-hand slope of the least concave majorant of the cumulative sum diagram defined by the set of points $\{(F_N(\overline x_k),\Lambda_N(\overline x_k)),k=0,\dots,K\}$ where $\overline x_0=0$. We define the corresponding inverse estimator as follows:
\begin{equation}\label{def: UN}
U_N(a)=\argmax_{u\in\{\overline x_0,\dots,\overline x_K\}}\{\Lambda_N(u)-aF_N(u)\}
\end{equation}
where $\overline x_0=0$, argmax denotes the greatest location of the maximum, and where we recall that for every nonincreasing left-continuous function $h:[0,1]\to\RR$, the generalized inverse of $h$  is defined as: for every $a\in\RR$, $h^{-1}(a)$ is the greatest  $t\in[0,1]$ that satisfies $h(t)\geq a$, with the convention that the supremum of an empty set is zero. To see that $U_N=\widehat\mu_N^{-1}$, note that from the characterization above of $\widehat\mu_N$ as the slope of a least concave majorant,  it follows  that for all $a\in\RR$ and $t\in(0,1]$, we have the equivalences
 \[
\begin{split}
\widehat\mu_N(t)<a
&\Leftrightarrow
\left(\exists \overline x_i<t\right)
\left(\forall \overline x_j\geq t\right):
\frac{\Lambda_N(\overline x_j)-\Lambda_N(\overline x_i)}{F_N(\overline x_j)-F_N(\overline x_i)}<a\\
&\Leftrightarrow
\left(\exists \overline x_i<t\right)
\left(\forall \overline x_j\geq t\right):
\Lambda_N(\overline x_j)-aF_N(\overline x_j)<\Lambda_N(\overline x_i)-aF_N(\overline x_i)\\
&\Leftrightarrow
\argmax_{u\in\{\overline x_0,\dots,\overline x_K\}}
\left\{
\Lambda_N(u)-aF_N(u)
\right\}
<t
\end{split}
\]
whereas for $t=0$, we have the equivalence
 \[
\begin{split}
\widehat\mu_N(0)< a
&\Leftrightarrow
\argmax_{u\in\{\overline x_0,\dots,\overline x_K\}}
\left\{
\Lambda_N(u)-aF_N(u)
\right\}=0.
\end{split}
\]
We study below the asymptotic properties of $U_N(a)$ for arbitrary $a$ and use these to deduce the asymptotic properties of $\widehat\mu_N(t)$ for a fixed $t\in(0,1)$ using the switch relation
\begin{equation}\label{eq: switch}
\widehat\mu_N(t)\geq a\Longleftrightarrow t\leq U_N(a),
\end{equation}
that holds for all $t\in(0,1]$ and $a\in\RR$. 
\newline
\newline
It will be useful to also record similar characterizations of the global estimator $\hat{\mu}_{N,G}$ of $\mu$, for the sake of completeness. Recall that the global estimator is the isotonic estimator that we would compute if all the data $\{X_i, Y_i\}_{i=1}^N$ could have been brought over (or were already there) on a central server. 
Letting $\Lambda_{N,G}(t) = N^{-1}\sum_{i=1}^N\,Y_i \mathds{1}_{X_i\leq t}$, for $a \in \mathbb{R}$, define 
\begin{equation}
\label{global-inv-est}
U_{N,G}(a) = \argmax_{u \in [0,1]}\{\Lambda_{N,G}(u)-aF_N(u)\} \,.
\end{equation}
Then $U_{N,G}(a) = \hat{\mu}_{N,G}^{-1}(a)$ and similar to the pooled estimator, we have the following characterization: 
\begin{equation}\label{eq:switchglobal}
\widehat\mu_{N,G}(t)\geq a\Longleftrightarrow t\leq U_{N,G}(a),
\end{equation}
that holds for all $t\in(0,1]$ and $a\in\RR$.
%
%Let $\Lambda_N$ be the piecewise-linear process on $[0,1]$ such that
%$$\Lambda_N\left(\sum_{k=1}^jw_k\right)=\sum_{k=1}^jw_k\overline y_{k}$$
%for all $j=1,\dots,K$, and $\Lambda_N(0)=0$.
%We then have
%\begin{equation}\label{brunk0}
%\widehat y_j=\widehat \mu_N(\overline x_{j})=\widehat\lambda_N\left(\sum_{k=1}^jw_k\right)=\widehat\lambda_N\left(F_N(\overline x_j)\right)
%\end{equation}
%for all $ j=1,\dots,K$, where $F_N$ is the empirical distribution function corresponding to $X_1,\dots,X_N$ and where $\widehat\lambda_N$ is the left-hand slope of the least concave majorant of $\Lambda_N.$
%We define the corresponding inverse estimator as follows:
%$$U_N(a)=\argmax_{u\in[0,1]}\{\Lambda_N(u)-au\}$$
%where argmax denotes the greatest location of the maximum, and where we recall that for every nonincreasing left-continuous function $h:[0,1]\to\RR$, the generalized inverse of $h$  is defined as: for every $a\in\RR$, $h^{-1}(a)$ is the greatest  $t\in[0,1]$ that satisfies $h(t)\geq a$, with the convention that the supremum of an empty set is zero. Then, $U_N=\lambda_N^{-1}$ and
%\begin{equation}\label{invmm}
%\widehat \mu_N^{-1}=F_N^{-1}\circ \widehat U_N.
%\end{equation}
\section{Asymptotic properties of the new estimators} 
In the sequel, we denote by $g$ the generalized inverse of $\mu$ and by $\EE^X$ the conditional expectation given $X_1,\dots,X_N$. Being the inverse  of $\mu$, $g$ is only defined on the interval $[\mu(1),\mu(0)]$. In the sequel, we expand $g$ to the whole real line by setting $g(a)=0$ for all $a>\mu(0)$ and $g(a)=1$ for all $a<\mu(1)$. 
\newline
Furthermore, for all $x\geq 0$, $[x]$ denotes the integer part of $x$. We denote by $F_X$ the mixing distribution function
\begin{equation}
\label{def: FX}
F_X(x)=\sum_{j=1}^m\frac{n_j}N F_{Xj}(x).
\end{equation}
Note that the function depends on $N$ but that for notational convenience, this is not made explicit in the notation. 
\newline
To develop the asymptotic properties of the proposed estimator, we will impose some further conditions on the model. These are:
\begin{itemize} 
\item [A1.] Assume that $F_X$ has a density function $f_X$ on $[0,1]$ that satisfies 
\begin{equation}\label{eq: HfX}
C_1< \inf_{t\in[0,1]}f_X(t)\leq \sup_{t\in[0,1]}f_X(t) \leq C_2
\end{equation}
for some positive numbers $C_1$ and $C_2$ that do not depend on $N$.

\item[A2.]  With $\varepsilon_i=Y_i-\mu(X_i)$ for all $i=1,\dots,N$, assume that there exists $\sigma>0$ such that $\EE[\varepsilon_i^2|X_i]\leq \sigma^2$ for all $i$, with probablity one. 

\item[A3.] The regression function $\mu$ satisfies: 
\begin{equation}\label{eq: Hmu}
C_3< \left\vert\frac{\mu(t)-\mu(x)}{t-x}\right\vert < C_4\mbox{ for all }t\neq x\in[0,1]\,,
\end{equation}
for positive numbers $C_3$ and $C_4$. 

\item[A4.] The number of bins $K$ satisfies $K^{-1}=o(N^{-1/3})$ and  there exists  $\lambda\in(0,1]$ that may depend on $N$ and satisfies 
\begin{equation}\label{eq: lambdaj}
 \min_{1\leq j\leq m}\frac {n_j}N\geq \lambda>0\mbox{ and } \liminf_{N\to\infty}N^{1/3}\lambda(\log N)^{-3}=\infty.
\end{equation} 
\end{itemize} 
\noindent
{\bf Remarks on the assumptions:} Assumption (A1) is fulfilled for instance if each $F_{Xj},\ j=1,\dots,m$ has a density function $f_{Xj}$ such that
$C_1 < f_{Xj}(x) < C_2\mbox{ for all }x\in[0,1]$. Note that Assumption (A3) is weaker than differentiability, it implies that $\mu$ is both Lipschitz and so to 
speak \emph{inverse Lipschitz}. It also implies that the inverse function $g$ defined above is continuous. Assumption (A4) is critical to recovering the Chernoff-type asymptotics for the pooled estimator; that $K$ grows faster than $N^{1/3}$ ensures that the data are averaged over bins of length smaller than $N^{-1/3}$, so that the isotonic algorithm operating on these averages at the central machine can still recover the $N^{-1/3}$ convergence rate. If $K$ were to grow exactly at the rate $N^{1/3}$ or slower, the pooled estimator would no longer demonstrate Chernoff-type cube-root asymptotics. Furthermore, in (A4), we assume that  the proportion $n_j/N$ of observations from the $j$-th sample is at least of order $N^{-1/3}(\log N)^{3}$. This also plays a critical role in the subsequent analysis. Since, $$1=\sum_{j=1}^m\frac{n_j}N\geq m\min_{1\leq j\leq m}\frac{n_j}N,$$ the conditions in \eqref{eq: lambdaj} imply that the number $m$ of different sub-samples cannot grow to fast: we must have $m\ll N^{1/3}(\log N)^{-3}.$
\subsection{Uniformly Bounded MSE Property of the New Estimators} 
%{\color{red}Some text will be added later on. I'll attend to this later. Also, we should add a parallel result as a second part of this theorem for the global estimator (without proof). In a remark, we should refer to BDS. We should note that Theorem 4.1 of BDS establishes the uniformly bounded maximal risk property of the global estimator under a strong (sub-Gaussian) assumption for the residual in the regression model, while this is not really necessary I think, since we don't need it for the inverse pooled estimator below. We should not need it for the inverse global estimator either, correct? \comC{no, we don't} The point is, the uniformly bounded maximal risk (MSE) property holds under much weaker conditions on the residuals.} 
%\newline
%\newline
{\bf The Inverse Problem:} We first demonstrate that the new estimator in the inverse problem exhibits \emph{uniformly bounded maximal risk (MSE)} over an appropriate class of models, as $N$ grows to $\infty$. This is an analogue of the first result in Theorem 4.1 of BDS for the \emph{global isotonic estimator of the inverse function}, though it is established here under weaker conditions. For this task, we denote by ${\cal F}_1$ the class of non-increasing  functions $\mu$ on $[0,1]$ that satisfy \eqref{eq: Hmu} and $\sup_t|\mu(t)|\leq C_5$, where $C_5>0$ is a positive number. 
\begin{theorem}\label{theo: supereffinv}
Under assumptions (A1) through (A4), there exists $C>0$ that depends only on $\sigma^2,C_1,C_2,C_3,C_4$ such that  for all $a\in\RR$,
\begin{equation*}
\limsup_{N\to\infty}\sup_{\mu\in{\cal F}_1}N^{2/3}\EE_\mu(U_N(a)-\mu^{-1}(a))^2\leq C.
\end{equation*}
\end{theorem}
The proof of the above theorem relies on a number of preliminary results which are presented, next. In the remainder of this section, we assume that assumptions (A1) to (A4) are always satisfied (though some results may require only a subset of these assumptions). Additional assumptions will be imposed 
when required. 
\begin{lemma}\label{lem: EXLambda}
%Assume that $F_X$ has a density function $f_X$ on $[0,1]$,  that $\mu$ is decreasing on $[0,1]$, and that there exist positive numbers $C_1,C_2,C_3,C_4$ such that
%\begin{equation}\label{eq: HfX}
%C_1< \inf_{t\in[0,1]}f_X(t)\leq \sup_{t\in[0,1]}f_X(t) \leq C_2
%\end{equation}
%and 
%%\begin{equation}\label{eq: Hmu}
%C_3< \left\vert\frac{\mu(t)-\mu(x)}{t-x}\right\vert < C_4\mbox{ for all }t\neq x\in[0,1].
%\end{equation}
%Assume moreover that $K^{-1}=o(N^{-1/3})$ and  there exists  $\lambda\in(0,1]$ (that may depend on $N$) such that
%\begin{equation}\label{eq: lambdaj}
% \min_{1\leq j\leq m}\frac {n_j}N\geq \lambda>0\mbox{ and } \liminf_{N\to\infty}N^{1/3}\lambda(\log N)^{-3}=\infty.
%\end{equation}
Let $\theta>0$ be arbitrary. Then,  
there exist (i) a number $c>0$ that depends only on $C_1,C_3$,  (ii) an integer $N_0>0$ that depends only on $C_1,C_2,C_3,C_4,\theta$ and (iii) an event ${\cal E}_N$ that depends only on $C_2$,  such that for all $N\geq N_0$, we have $\PP({\cal E}_N)\geq 1- N^{-\theta}$  and on ${\cal E}_N$,
\begin{equation}\label{eq: 5.1}
E^X\Lambda_N(u)-E^X\Lambda_N\left(\frac{[Kg(a)]}{K}\right)-a\left(F_N(u)-F_N\left(\frac{[Kg(a)]}{K}\right)\right)\leq-c(u-g(a))^2
\end{equation}
for all $a\in\RR$ and all $u\in\{\overline x_0,\dots,\overline x_K\}$ such that $|u-g(a)|\geq N^{-1/3}$.
\end{lemma}
The proof of this lemma is long and technical and is available in the appendix. The next result gives a polynomial tail bound on the estimation error $U_N(a) - g(a)$ over a high-probability set that is eventually used to bound the MSE. 
%{\color{red}conditions on $\mu$ missing, need to be inserted}
%{\color{red}This is a minor point but it seems to me that we are also implicitly using the fact that $o((g(a) - u)^2)/(g(a)-u)^2$ is also uniformly bounded for all $|g(a) - u| \geq N^{-1/3}$ with $u$ varying between 0 and 1 and $|g(a) - u| \geq N^{-1/3}$. \comC{I have clarified the point in blue above.} Because, from the above display, for $\tilde{c} < C_1/2$, we have: 
%\[f(a,u) \geq \tilde{c}(g(a) - u)^2 + (g(a)-u)^2[ (C_1/2 - \tilde{c}) + o((g(a)-u)^2)/(g(a)-u)^2)]\,,\]
%and we can say that the $f(a,u) \geq \tilde{c}(g(a)-u)^2$ only if the second term on the right-side is positive for all $u$ and $|g(a)-u|\geq N^{-1/3}$. This would mean that we need the term $((C_1/2 - \tilde{c}) + o((g(a)-u)^2)/(g(a)-u)^2$ to be always non-negative; because $u$ can vary between $\overline{x}_0$ and $\overline{x}_K$, i.e. from close to 0 to close to 1 and $g(a)$ is a fixed point between 0 and 1, it is not clear what this term is doing when $g(a)$ and $u$ are not close. Unless, by uniformity you imply that 
%$$\lim_{N \rightarrow \infty} \sup_{u, a: |g(a)-u| \geq N^{-1/3}} o((g(a)-u)^2)/(g(a)-u)^2 = 0 \,.$$ 
%On going again through the proof steps, it seems to me that this is exactly what is happening. So I think I was initially misinterpreting the meaning of the $o((g(a)-u)^2$ terms. I think it might be helpful to just make the meaning of the uniformity explicit via the above equation.} 
\begin{lemma}\label{lem: tail boundUN}
With ${\cal E}_N$ and $N_0$ taken from Lemma \ref{lem: EXLambda}, there exists $C>0$ that depends only on $\sigma^2,C_1,C_2,C_3,C_4$ such that  for all $a\in\RR$ and $x>0$,
\begin{equation}\label{eq: tail boundUN}
 \PP\left(|U_N(a)-g(a)|\geq x,{\cal E}_N\right)\leq \frac{C}{Nx^3}
\end{equation}
for all $N\geq N_0$.
\end{lemma}
\noindent
{\bf Proof.} The inequality in the lemma is obvious for $x\in(0,N^{-1/3})$ since for such $x$'s, it suffices to choose $C\geq 1$ so that the right hand side is larger than one. Hence, in the sequel we consider $x\geq N^{-1/3}$. For all $a\in\RR$ and all $u\in\{\overline x_0,\dots,\overline x_K\}$ such that $|u-g(a)|\geq x$, define $e(a,u)$ as in
\eqref{eq: eau} and $M_N(u)=\Lambda_N(u)-E^X(\Lambda_N(u)).$ The characterization in  \eqref{def: UN} proves the following inclusion of events:
\begin{eqnarray*}
&&\left\{U_N(a)-g(a)\geq x\right\}\\
&&\qquad \subset\left\{\max_{u\in\{\overline x_0,\dots,\overline x_K\},\ u-g(a)\geq x}\{\Lambda_N(u)-aF_N(u)\}\geq \Lambda_N\left(\frac{[Kg(a)]}{K}\right)-aF_N\left(\frac{[Kg(a)]}{K}\right)\right\}\\
&&\qquad=\left\{\max_{u\in\{\overline x_0,\dots,\overline x_K\},\ u-g(a)\geq x}\left\{M_N(u)- M_N\left(\frac{[Kg(a)]}{K}\right)+e(a,u)\right\}\geq 0\right\}.
\end{eqnarray*}
Since $x\geq N^{-1/3}$, combining this with Lemma \ref{lem: EXLambda} shows that there exists $c>0$ that depends only on $C_1,C_3$ such that
\begin{eqnarray*}
&&\PP\left(U_N(a)-g(a)\geq x,{\cal E}_N\right)\\
&&\qquad\leq\PP\left(\max_{u\in\{\overline x_0,\dots,\overline x_K\},\ u-g(a)\geq x}\left\{M_N(u)- M_N\left(\frac{[Kg(a)]}{K}\right)-c(u-g(a))^2\right\}\geq 0\right)
\end{eqnarray*}
for $N\geq N_0$.
Hence,
\begin{eqnarray}\notag \label{eq: sumk}
&&\PP\left(U_N(a)-g(a)\geq x,{\cal E}_N\right)\\ \notag
&&\qquad\leq\sum_{k\geq 0}\PP\left(\max_{u\in\{\overline x_0,\dots,\overline x_K\},\ u-g(a) \in[x2^k,x2^{k+1}]}\left\{M_N(u)- M_N\left(\frac{[Kg(a)]}{K}\right)-c(u-g(a))^2\right\}\geq 0\right)\\ 
&&\qquad\leq\sum_{k\geq 0}\PP\left(\max_{u\in\{\overline x_0,\dots,\overline x_K\},\ u-g(a) \in[0,x2^{k+1}]}\left\{M_N(u)- M_N\left(\frac{[Kg(a)]}{K}\right)\right\}\geq c(x2^k)^2\right).
\end{eqnarray}
Let $\PP^X$ denote the conditional probability given $X_1,\dots,X_N$. By definition, for all $u\in\{\overline x_0,\dots,\overline x_K\}$ we have
\begin{equation}\label{eq: defMN}
M_N(u)=\frac 1N\sum_{i=1}^N\varepsilon_i\mathds{1}_{X_i\leq u}
\end{equation}
where $\varepsilon_i=Y_i-\mu(X_i)$. The process $M_n$ can be extended to all $u\in\RR$ using the same definition as above. Then,  $M_N$ is a centered martingale under $\PP^X$ that satisfies
\begin{equation}\label{eq: varMN}
\EE^X\left(M_N(u)-M_N(v)\right)^2 =\frac 1{N^2}\sum_{i=1}^N\EE^X(\varepsilon_i^2)\mathds{1}_{u<X_i\leq v}\leq \frac {\sigma^2}{N}(F_N(u)-F_N(v))
\end{equation}
for all $u\leq v$, using that $\EE^X(\varepsilon_i^2)\leq\sigma^2$ for all $i$ by assumption. Hence, it follows from the Doob inequality that for all $k\geq 0$,
\begin{eqnarray*}
&&\PP^X\left(\max_{u\in\{\overline x_0,\dots,\overline x_K\},\ u-g(a) \in[0,x2^{k+1}]}\left\{M_N(u)- M_N\left(\frac{[Kg(a)]}{K}\right)\right\}\geq c(x2^k)^2\right)\\
&&\qquad\leq \sigma^2\frac{F_N\left(g(a)+x2^{k+1}\right)-F_N\left(\frac{[Kg(a)]}{K}\right)}{c^2N(x2^k)^4}.
\end{eqnarray*}
Taking the expectation on both sides of the preceding inequality yields for large enough $N$ that
\begin{eqnarray*}
&&\PP\left(\max_{u\in\{\overline x_0,\dots,\overline x_K\},\ u-g(a) \in[0,x2^{k+1}]}\left\{M_N(u)- M_N\left(\frac{[Kg(a)]}{K}\right)\right\}\geq c(x2^k)^2\right)\\
&&\qquad\leq \frac{\sigma^2\left\{F_X\left(g(a)+x2^{k+1}\right)-F_X\left(\frac{[Kg(a)]}{K}\right)\right\}}{c^2N(x2^k)^4}\\
&&\qquad\leq \frac{\sigma^2 C_2(x2^{k+1}+K^{-1})}{c^2N(x2^k)^4}\\
&&\qquad\leq \frac{2\sigma^2 C_2x2^{k+1}}{c^2N(x2^k)^4},
\end{eqnarray*}
where $C_2$ is taken from \eqref{eq: HfX}. For the penultimate inequality, we used that  $x2^{k+1}\geq N^{-1/3}$ for all $k$ whereas $K^{-1}=o(N^{-1/3})$, implying that $K^{-1}\leq x2^{k+1}$ for all $k$ provided that $N$ is sufficiently large. Putting the previous inequality in \eqref{eq: sumk} we obtain that for sufficiently large $N$,
\begin{eqnarray*}
\PP\left(U_N(a)-g(a)\geq x,{\cal E}_N\right)\leq\sum_{k\geq 0}\frac{4\sigma^2C_2}{c^2N(x2^k)^3}.
\end{eqnarray*}
Since $C:=\sum_{k\geq 0}2^{-3k}$ is finite, we conclude that 
\begin{eqnarray*}
\PP\left(U_N(a)-g(a)\geq x,{\cal E}_N\right)\leq\frac{4\sigma^2C_2C}{c^2Nx^3}.
\end{eqnarray*}
Similar arguments show that
\begin{eqnarray*}
\PP\left(g(a)-U_N(a)\geq x,{\cal E}_N\right)\leq\frac{4\sigma^2C_2C}{c^2Nx^3},
\end{eqnarray*}
and therefore,
\begin{eqnarray*}
\PP\left(|g(a)-U_N(a)|\geq x,{\cal E}_N\right)\leq\frac{8\sigma^2C_2C}{c^2Nx^3}.
\end{eqnarray*}
The lemma follows.\hfill{$\Box$}
\newline
\newline
We are now ready to prove the theorem. 
\newline
\newline
{\bf Proof of Theorem \ref{theo: supereffinv}.} Fix $\mu\in{\cal F}_1$ arbitrarily. Since both $U_N$ and $\mu^{-1}$ take values in $[0,1]$, we have $|U_N(a)-\mu^{-1}(a)|\leq 1$ for all $a$ and therefore, with $\overline{\cal E}_N$ the complementary event to ${\cal E}_N$ taken from Lemma \ref{lem: EXLambda}, where we set $\theta=2/3$, we have
\begin{eqnarray}\label{eq: se}\notag
\EE_\mu\left(|U_N(a)-\mu^{-1}(a)|^2\mathds{1}_{\overline{\cal E}_N}\right)
&\leq&\PP_\mu\left(\overline{\cal E}_N\right)\\
&\leq &N^{-2/3}
\end{eqnarray}
for $N$ sufficiently large. On the other hand, 
it follows from the Fubini theorem that
\begin{eqnarray*}
\EE_\mu\left(|U_N(a)-\mu^{-1}(a)|^2\mathds{1}_{{\cal E}_N}\right)
&=&\int_0^\infty\PP_\mu\left(|U_N(a)-\mu^{-1}(a)|>\sqrt x,{\cal E}_N\right)dx\\
&=&\int_0^\infty2y\PP_\mu\left(|U_N(a)-\mu^{-1}(a)|>y,{\cal E}_N\right)dy\\
&\leq &\int_0^\infty2y\left(\frac{C}{Ny^3}\wedge 1\right)dy.
\end{eqnarray*}
For the last inequality, we used \eqref{eq: tail boundUN} together with the fact that a probability cannot be larger than one. Hence,
\begin{eqnarray*}
\EE_\mu\left(|U_N(a)-\mu^{-1}(a)|^2\mathds{1}_{{\cal E}_N}\right)
&\leq &\int_0^{N^{-1/3}}2ydy+\int_{N^{-1/3}}^\infty\frac{2C}{Ny^2}dy\\
&\leq& N^{-2/3}\left(1+2C\right).
\end{eqnarray*}
Combining with \eqref{eq: se} yields
\begin{eqnarray*}
\EE_\mu\left(|U_N(a)-\mu^{-1}(a)|^2\right)
&\leq& N^{-2/3}\left(2+2C\right),
\end{eqnarray*}
which completes the proof of the Theorem (by taking $C$ to be $2 + 2C$) where $C$ is the constant from Lemma \eqref{lem: tail boundUN}.
\hfill{$\Box$}
\newline
\newline 
{\bf The Direct Problem:}  We now establish an analogue of the second result in Theorem 4.1 of BDS to demonstrate that the new estimator fixes the super-effciency phenomenon in the direct problem as well, i.e. it has bounded uniform MSE as $N \rightarrow \infty$ over the class $\mathcal{F}_1$. 
\newline
Denote by $\widetilde F_N$ the step function on $[0,1]$ such that $\widetilde F_N(\overline x_k)= F_N(\overline x_k)$ for all $k=0,\dots,K$, and $\widetilde F_N$ is constant on all intervals $[\overline x_{k-1},\overline x_k)$ for $k=1,\dots,K$. We denote by $\widetilde F_N^{-1}$ the corresponding inverse function:
\begin{equation*}
\widetilde F_N^{-1}(t)=\inf\{x\in[0,1]\mbox{ such that }\widetilde F_N(x)\geq t\}.
\end{equation*}
Since $\widetilde F_N^{-1}\circ \widetilde F_N(\overline x_k)=\overline x_k$ for all $k=0,\dots,K$, it follows from the characterization in  \eqref{def: UN} that
\begin{equation}\label{eq: VtoU}
U_N(a)=\widetilde F_N^{-1}(V_N(a))
\end{equation}
for all $a\in\RR$, where
\begin{equation*}
V_N(a)=\argmax_{u\in\{\widetilde F_N(\overline x_0),\dots,\widetilde F_N(\overline x_K)\}}\{\Lambda_N\circ \widetilde F_N^{-1}(u)-au\}.
\end{equation*}
The following lemma provides tail bound probabilities for $V_N$.
\begin{lemma}\label{lem: tail boundVN}
With $\varepsilon_i=Y_i-\mu(X_i)$ for all $i=1,\dots,N$, assume that there exists $\sigma>0$ such that $\EE[\varepsilon_i^p|X_i]\leq \sigma^p$ for all $i$ and some $p\geq 2$, with probablity one. Assume that $F_X$ has a density function $f_X$ on $[0,1]$ that satisfies \eqref{eq: HfX} for some positive numbers $C_1,C_2$.  Then, there exists  $C>0$ that depends only on $p,C_2$ and $\sigma$ such that
\begin{eqnarray*}
\PP\left(V_N(a)\geq x\right)\leq\frac{C}{N^{p/2}x^{p-1}(a-\mu(0))^p}
\end{eqnarray*}
for all $a>\mu(0)$ and
\begin{eqnarray*}
\PP\left(1-V_N(a)\geq x\right)\leq\frac{C}{N^{p/2}x^{p-1}(\mu(1)-a)^p}.
\end{eqnarray*}
for all $a<\mu(1)$.
\end{lemma}
For a proof of this lemma, see the Appendix. 
\newline
\newline
The following theorem establishes the desired property for the pooled direct estimator.  
\begin{theorem}\label{theo: supereffdirect}
Fix $\delta\in(0,1/2)$. Then, there exists $C>0$ that depends only on $\sigma,p,C_1,C_2,C_3,C_4,C_5,\delta$ such that  for all $t\in[\delta,1-\delta]$
\begin{equation*}
\limsup_{N\to\infty}\sup_{\mu\in{\cal F}_1}N^{2/3}\EE_\mu(\widehat\mu_N(t)-\mu(t))^2\leq C.
\end{equation*}
\end{theorem}
\medskip
%{\color{red}Similar to the inverse function case, we should add a second part about the bounded maximal risk of the global estimator under similar assumptions.} 
%\newline
{\bf Proof.} Similar to the proof of Theorem \ref{theo: supereffinv} for the inverse problem, 
we would like to restrict ourselves to the event ${\cal E}_N$ from Lemma \ref{lem: EXLambda}, where $\theta$ can be chosen arbitrarily large.
However, we do not have an analogue of \eqref{eq: se} for the direct problem since $\widehat \mu_N$ is not bounded as is $U_N$. Hence, we first prove that $\widehat\mu_N$ remains bounded by a power of $N$ apart possibly on a negligible set. For this task, consider an arbitrary $A>0$ such that $A+\mu(0)>0$, and note that for all $t\in[0,1]$, and all non-increasing functions $\mu$ on $[0,1]$, we have
\begin{eqnarray*}
\EE_\mu\left[\widehat\mu_N^2(t)\mathds{1}_{\widehat\mu_N(t)>A+\mu(0)}\right]\leq
\EE_\mu\left[\widehat\mu_N^2(0)\mathds{1}_{\widehat\mu_N(0)>A+\mu(0)}\right].
\end{eqnarray*}
Hence, it follows from the Fubini theorem that for all non-increasing $\mu\in{\cal F}_1$,
\begin{eqnarray*}
\EE_\mu\left[\widehat\mu_N^2(t)\mathds{1}_{\widehat\mu_N(t)>A+\mu(0)}\right]
&\leq&
\int_{0}^\infty\PP_\mu(\widehat\mu_n(0)\mathds{1}_{\widehat\mu_N(0)>A+\mu(0)}>\sqrt y)dy\\
&=& (A+\mu(0))^2\PP_\mu(\widehat\mu_N(0)>A+\mu(0))+
\int_{A+\mu(0)}^\infty2y\PP_\mu(\widehat\mu_N(0)>y)dy.
\end{eqnarray*}
Note that if $\widehat\mu_N(0)>y$ for some $y\in\RR$, then for the inverse we must have $U_N(y)>0$. Since $U_N$ can only assume values in the set of jump points of $\widehat \mu_N$ it is of the form $\overline{x}_k = k/K$ for some $k \geq 1$. Next, $\widehat\mu_N$ can have jumps only at those $\overline x_k$ where $\widetilde F_N$ has a jump, i.e. $\widetilde F_N(\overline x_k)>\widetilde F_N(\overline x_{k-1})$. Since the size of a jump of $\widetilde F_N$ is at least $N^{-1}$, we must have $\widetilde F_N(\overline x_k)\geq N^{-1}$ and therefore, $F_N(\overline x_k)=\widetilde F_N(\overline x_k)\geq N^{-1}$. Thus,  
$$V_N(y)=\widetilde F_N(U_N(y))=F_N(U_N(y)) = F_N(\overline{x}_k)\geq N^{-1},$$
\newline
implying that for all $\mu\in{\cal F}_1$,
\begin{eqnarray*}
\EE_\mu\left[\widehat\mu_N^2(t)\mathds{1}_{\widehat\mu_N(t)>A+\mu(0)}\right]
&\leq&
(A+\mu(0))^2\PP_\mu(V_N(A+\mu(0))\geq N^{-1})+
\int_{A+\mu(0)}^\infty2y\PP_\mu(V_N(y)\geq N^{-1})dy.
\end{eqnarray*}
With $C$ taken from Lemma \eqref{lem: tail boundVN} where it is assumed that $p>2$, we arrive at
\begin{eqnarray*}
\EE_\mu\left[\widehat\mu_N^2(t)\mathds{1}_{\widehat\mu_N(t)>A+\mu(0)}\right]
&\leq&
CN^{-1+p/2}(A+\mu(0))^2A^{-p}+ 2CN^{-1+p/2}
\int_{A+\mu(0)}^\infty y(y-\mu(0))^{-p}dy\\
&=&
CN^{-1+p/2}(A+\mu(0))^2A^{-p}+2CN^{-1+p/2}\left\{\frac{A^{2-p}}{p-2}+\mu(0)\frac{A^{1-p}}{p-1}\right\}\\
&\leq&
CN^{-1+p/2}(A+C_5)^2A^{-p}+2CN^{-1+p/2}\left\{\frac{A^{2-p}}{p-2}+C_5\frac{A^{1-p}}{p-1}\right\}.
\end{eqnarray*}
With $A=N^{(3p-2)/(6(p-2))}$, this proves that there exists $C'>0$ that depends only on $\sigma$, $p$, $C_2$ and $C_5$ such that 
%({\color{red} we don't want 
%the bound to depend on a particular $\mu(0)$, but the bound really depends on the maximal absolute value of $\mu(0)$ which is finite over the class $\mathcal{F}_1$, so this is what we should probably say} \comC{You're right. This is fixed by the last line in green in the previous display.}) such that
\begin{eqnarray*}
\EE_\mu\left[(\widehat\mu_N(t))^2\mathds{1}_{\widehat\mu_N(t)>A+\mu(0)}\right]
&\leq&C'N^{-2/3}
\end{eqnarray*}
for all $t\in[0,1]$ and $\mu\in{\cal F}_1$. Now, with $A=N^{(3p-2)/(6(p-2))}$,
\begin{eqnarray}\notag
\EE_\mu\left[(\widehat\mu_N(t)-\mu(t))^2\mathds{1}_{\widehat\mu_N(t)>A+\mu(0)}\right]
&\leq&\EE_\mu\left[2\left(\widehat\mu^2_N(t)+\mu^2(t)\right)\mathds{1}
_{\widehat\mu_N(t)>A+\mu(0)}\right]\\ \notag
&\leq&2C'N^{-2/3}+2\max\{|\mu(0)|,|\mu(1)|\}^2\PP\left(\widehat\mu_N(0)>A+\mu(0)\right)\\ \notag
&\leq&2C'N^{-2/3}+2\max\{|\mu(0)|,|\mu(1)|\}^2\PP\left(V_N(A+\mu(0))\geq N^{-1}\right),
\end{eqnarray}
similar as above, whence
\begin{eqnarray}\label{eq: A+mu0}\notag
\EE_\mu\left[(\widehat\mu_N(t)-\mu(t))^2\mathds{1}_{\widehat\mu_N(t)>A+\mu(0)}\right]
&\leq&2C'N^{-2/3}+2C{C_5}^2N^{-1+p/2}A^{-p}\\
&\leq&C''N^{-2/3}
\end{eqnarray}where $C''$ depends only on $\sigma,p,C_2,$ and $C_5$. 
This enables us to restrict to the event ${\cal E}_N$ of Lemma \ref{lem: EXLambda}, provided that $\theta$ is chosen sufficiently large in the lemma. Indeed, with $\theta>(5p-6)/(3(p-2))$, the previous inequality implies that with $A=N^{(3p-2)/(6(p-2))}$ and $N$ sufficiently large,
\begin{eqnarray}\notag
\EE_\mu\left[\{(\widehat\mu_N(t)-\mu(t))_+\}^2\mathds{1}_{\overline{\cal E}_N}\right]
&\leq&(A+\mu(0)-\mu(t))^2\PP_\mu(\overline{\cal E}_N)+C''N^{-2/3}\\
&\leq&(A+2C_5)^2\PP_\mu(\overline{\cal E}_N)+C''N^{-2/3}\\
&\leq& 2C''N^{-2/3}
\end{eqnarray}
%({\color{red}I get the same result as in the above equations but in a slightly different way. See Distributed Computing.pdf in the shared folder. \comC{Is it clearler with the additional line in the previous display?}Also, while it is clear from the context, at a first glance it seems you are taking the positive part of the squared difference, whereas it is really the square of the positive part. Maybe we could change the notation slightly.} \comC{I have added brakets to avoid ambiguity}) 
for  all $t\in[0,1]$  and all $\mu\in{\cal F}_1$. It can be shown similarly that for $N$ sufficiently large,
\begin{eqnarray}\notag
\EE_\mu\left[\{(\widehat\mu_N(t)-\mu(t))_-\}^2\mathds{1}_{\overline{\cal E}_N}\right]
&\leq& 2C''N^{-2/3}
\end{eqnarray}
for  all $t\in[0,1]$  and all $\mu\in{\cal F}_1$, implying that
\begin{eqnarray}\notag
\limsup_{N\to\infty}\sup_{\mu\in{\cal F}_1}N^{2/3}\EE_\mu\left[(\widehat\mu_N(t)-\mu(t))^2\mathds{1}_{\overline{\cal E}_N}\right]
&\leq& 4C''
\end{eqnarray}
for  all $t\in[0,1]$. Hence, it now suffices to prove that there exists $C>0$ that depends only on $\sigma,p,C_1,C_2,C_3,C_4,C_5,\delta$ such that
\begin{equation}\label{eq: seEN}
\limsup_{N\to\infty}\sup_{\mu\in{\cal F}_1}N^{2/3}\EE_\mu\left[(\widehat\mu_N(t)-\mu(t))^2\mathds{1}_{{\cal E}_N}\right]\leq C.
\end{equation}
To prove this, fix $\mu\in{\cal F}_1$ arbitrarily, and invoke the Fubini Theorem to obtain that
\begin{eqnarray}\notag\label{eq: Fubes}
\EE_\mu\left[\{(\widehat\mu_N(t)-\mu(t))_+\}^2\mathds{1}_{{\cal E}_N}\right]
&=&\int_0^\infty2y\PP_\mu\left(\widehat \mu_N(t)-\mu(t)\geq y,{\cal E}_N
\right)dy\\
&=&\int_0^\infty2y\PP_\mu\left(U_N(\mu(t)+y)\geq t,{\cal E}_N
\right)dy,
\end{eqnarray}
 using the switch relation \eqref{eq: switch} for the last equality. We split the above integral into the sum of two integrals and first consider
\begin{eqnarray*}
I_1=\int_0^{\mu(0)-\mu(t)}2y\PP_\mu\left(U_N(\mu(t)+y)\geq t,{\cal E}_N\right)dy.
\end{eqnarray*}
With $C_4$ taken from the definition of ${\cal F}_1$ we have
\begin{equation*}
t= \mu^{-1}(\mu(t))\geq \mu^{-1}(\mu(t)+y)+yC_4^{-1}
\end{equation*}
for all $t\in[0,1]$ and $y\in[0,\mu(0)-\mu(t)]$. Combining Lemma \ref{lem: tail boundUN} with the fact that a probability cannot be larger than one then yields 
\begin{eqnarray}\notag\label{eq: majI1}
I_1&\leq &N^{-2/3}+\int_{N^{-1/3}}^{\mu(0)-\mu(t)}2y\PP_\mu\left(U_N(\mu(t)+y)- \mu^{-1}(\mu(t)+y) \geq yC_4^{-1},{\cal E}_N\right)dy\\ \notag
&\leq &N^{-2/3}+\int_{N^{-1/3}}^{\infty}\frac{2CC_4^3}{Ny^2}dy\\
&\leq &N^{-2/3}\left(1+2CC_4^3\right).
\end{eqnarray}
Next, Lemma \ref{lem: tail boundUN} yields
\begin{eqnarray*}
I_2&:=&\int_{\mu(0)-\mu(t)}^\infty2y\PP_\mu\left(U_N(\mu(t)+y)\geq t,{\cal E}_N\right)dy\\
&\leq &\int_{\mu(0)-\mu(t)}^{\mu(0)-\mu(t)+N^{1/6}}2y\PP_\mu\left(U_N(0)\geq t,{\cal E}_N\right)dy+\int_{\mu(0)-\mu(t)+N^{1/6}}^\infty2y\PP_\mu\left(U_N(\mu(t)+y)\geq t,{\cal E}_N\right)dy
\\
&\leq &\frac{C}{Nt^3}\int_{\mu(0)-\mu(t)}^{\mu(0)-\mu(t)+N^{1/6}}2ydy+\int_{\mu(0)-\mu(t)+N^{1/6}}^\infty2y\PP_\mu\left(U_N(\mu(t)+y)\geq t,{\cal E}_N\right)dy,
\\
%&\leq &{\color{blue}\frac{2C}{Nt^3}(\mu(0)-\mu(t)+N^{1/6})N^{1/6}}+\int_{\mu(0)-\mu(t)+N^{1/6}}^\infty2y\PP_\mu\left(U_N(\mu(t)+y)\geq t,{\cal E}_N\right)dy\\
%&\leq &{\color{blue}\frac{3C}{\delta^3}N^{-2/3}} +\int_{\mu(0)-\mu(t)+N^{1/6}}^\infty2y\PP_\mu\left(U_N(\mu(t)+y)\geq t,{\cal E}_N\right)dy,
\end{eqnarray*}
where the first term on the right-hand side is equal to
\begin{eqnarray*}
\frac{C}{Nt^3}\left((\mu(0)-\mu(t)+N^{1/6})^2-(\mu(0)-\mu(t))^2\right)&=&
\frac{C}{Nt^3}\left(2(\mu(0)-\mu(t))N^{1/6}+N^{1/3}\right)\\
&\leq&\frac{C}{Nt^3}(4C_5N^{-1/6}+1)N^{1/3}\\
&\leq&\frac{2C}{\delta^3}N^{-2/3}
\end{eqnarray*}
for sufficiently large $N$, for all $t\geq\delta$ and $\mu\in{\cal F}_1$. Using the connection \eqref{eq: VtoU} between $U_N$ and $V_N$ yields
\begin{eqnarray*}
I_2&\leq &\frac{2C}{\delta^3}N^{-2/3}+\int_{\mu(0)-\mu(t)+N^{1/6}}^\infty2y\PP_\mu\left(V_N(\mu(t)+y)\geq \widetilde F_N(t),{\cal E}_N\right)dy,
\end{eqnarray*}
where $\widetilde F_N(t)=\widetilde F_N([Kt]K^{-1})=F_N([Kt]K^{-1})$ by definition of $\widetilde F_N$ and $F_N$. Regarding the proof of Lemma \ref{lem: EXLambda}, it can be seen that on ${\cal E}_N$ we have
$$\sup_{t\in[0,1]}\vert F_{N}(t)-F_X(t)\vert\leq C_2N^{-1/3}$$
whence
\begin{eqnarray*}
I_2&\leq &\frac{2C}{\delta^3}N^{-2/3}+\int_{\mu(0)-\mu(t)+N^{1/6}}^\infty2y\PP_\mu\left(V_N(\mu(t)+y)\geq F_X([Kt]K^{-1})-C_2N^{-1/3}\right)dy\\
&\leq &\frac{2C}{\delta^3}N^{-2/3}+\int_{\mu(0)-\mu(t)+N^{1/6}}^\infty2y\PP_\mu\left(V_N(\mu(t)+y)\geq C_1(t-K^{-1})-C_2N^{-1/3}\right)dy\\
&\leq &\frac{2C}{\delta^3}N^{-2/3}+\int_{\mu(0)-\mu(t)+N^{1/6}}^\infty2y\PP_\mu\left(V_N(\mu(t)+y)\geq C_1\delta/2\right)dy,
\end{eqnarray*}
for all $t\in[\delta,1-\delta]$, provided that $N$ is sufficiently large. Hence, it follows from Lemma \ref{lem: tail boundVN} where it is assumed that $p>2$, that
\begin{eqnarray*}
I_2
&\leq &\frac{2C}{\delta^3}N^{-2/3}+\frac{2^pC}{(C_1\delta)^{p-1}}\int_{\mu(0)-\mu(t)+N^{1/6}}^\infty
\frac{yN^{-p/2}}{(y+\mu(t)-\mu(0))^p}dy.
\end{eqnarray*}
For the integral on the right-hand side we have
\begin{eqnarray*}
&&\int_{\mu(0)-\mu(t)+N^{1/6}}^\infty
\frac{yN^{-p/2}}{(y+\mu(t)-\mu(0))^p}dy\\
&&\quad =\int_{\mu(0)-\mu(t)+N^{1/6}}^\infty
\frac{N^{-p/2}}{(y+\mu(t)-\mu(0))^{p-1}}dy+\int_{\mu(0)-\mu(t)+N^{1/6}}^\infty
\frac{(\mu(0)-\mu(t))N^{-p/2}}{(y+\mu(t)-\mu(0))^p}dy\\
&&\quad =\int_{N^{1/6}}^\infty
\frac{N^{-p/2}}{u^{p-1}}du+\int_{N^{1/6}}^\infty
\frac{(\mu(0)-\mu(t))N^{-p/2}}{u^p}du\\
&&\qquad\leq \frac{1}{p-2}N^{(1-2p)/3}+\frac{2C_5}{p-1}N^{(1-4p)/6}.
\end{eqnarray*}
Hence, we can find $\tilde C$ that depends only on $p,\sigma,C_1-C_5$ such that
$$I_2\leq \tilde CN^{-2/3}$$
for sufficiently large $N$, for all $\mu\in{\cal F}_1$ and $t\in[\delta,1-\delta]$.
\newline
Combining this with \eqref{eq: majI1} and \eqref{eq: Fubes} proves that there exists $C>0$ that depends only on $\sigma,p,C_1,C_2,C_3,C_4,C_5,\delta$ such that
\begin{equation}\notag
\limsup_{N\to\infty}\sup_{\mu\in{\cal F}_1}N^{2/3}\EE_\mu\left[(\widehat\mu_N(t)-\mu(t))_+^2\mathds{1}_{{\cal E}_N}\right]\leq C.
\end{equation}
It can be proved similarly that
\begin{equation}\notag
\limsup_{N\to\infty}\sup_{\mu\in{\cal F}_1}N^{2/3}\EE_\mu\left[(\widehat\mu_N(t)-\mu(t))_-^2\mathds{1}_{{\cal E}_N}\right]\leq C,
\end{equation}
which
completes the proof \eqref{eq: seEN}, and hence the proof of the theorem. \hfill{$\Box$}
\newline
\newline
In the next section, we show that under a fixed $\mu$, the new estimator recovers the asymptotic distribution of the global estimator with the same convergence rate.
\subsection{Asymptotic distributions} 
To establish asymptotic distributions for our new estimators, we make \emph{additional assumptions} in the case that the number $m$ of different samples goes to infinity, and we clarify the asymptotic setting further. 
\newline
\newline
When considering the case where $m$ is allowed to grow to infinity as $N\to\infty$, we assume that there is a sequence of unknown distinct distributions $\{P_j\}_{j\geq 1}$ such that our set of observations is part of an infinite sequence of pairs $\{(X_i,Y_i)\}_{i\geq 1}$, where for all $i$ the distribution of $(X_i,Y_i)$ takes the form $P_j$ for some $j\geq 1$. Hence, $m=m_N$ is the number of different distributions that appear across the first $N$ observations $(X_1,Y_1),\dots,(X_N,Y_N)$. To fix ideas, possibly rearranging the probabilities in the sequence $\{P_j\}_{j\geq 1}$, we assume without loss of generality in the sequel that for all $N$, the $m=m_N$ distributions 
that appear across the first $N$ observations are $P_1,\dots,P_m$. Note that the setting does not exclude that $m_N=1$ for all $N$, i.e. that all observations are drawn from the same distribution $P_1$. In the case where $m_N>1$ for sufficiently large $N$, it is not excluded that $m_N$ remains bounded.  In the sequel, for all $j\geq 1$, we denote by $\sigma_j$ the function such that
$$\sigma_j^2(u)=\EE[(Y-\mu(X))^2|X=u]$$
for all $u\in[0,1]$ and by $f_j$ the density function of $X$, which is assumed to exist, where $(X,Y)$ has distribution $P_j$. Then, the distribution function $F_X$ in \eqref{def: FX} has a density function $f_X$ on $[0,1]$ given by
\begin{eqnarray}\label{eq: fX}
f_X(u)=\sum_{j=1}^m\frac{n_j}N f_{j}(u)
\end{eqnarray}
for all $u\in[0,1]$. 
\newline
\newline
We next make the following technical assumptions.
\begin{itemize}

\item[$\tilde{A}_0$.] The functions $\{f_j\}$ are uniformly bounded in $j$ on the interval $[0,1]$. 

\item[$\tilde{A}_1$.] Let 
\begin{eqnarray*}
\omega(\delta)=\sup_{j\geq 1}\max\{\sup_{|u-v|\leq\delta}|\sigma_j^2(u)-\sigma_j^2(v)|,\sup_{|u-v|\leq\delta}|f_j(u)-f_j(v)|\}
\end{eqnarray*}
for all $\delta\geq 0$. Then, $\omega(\delta)\to 0$ as $\delta\to0$. 

\item[$\tilde{A}_2$.] The density function $f_X$ converges pointwise [and hence, uniformly] on $[0,1]$ as $N\to\infty$ to a continuous function $f_\infty$ that is bounded away from zero. This implies that \eqref{eq: HfX} holds for some positive numbers $C_1,C_2$ that do not depend on $N$, provided that $N$ is sufficiently large. 

\item[$\tilde{A}_3$.] The function $\sigma^2_X$ defined by $$\sigma^2_X(u):=\sum_{j=1}^m\frac{n_j}N \sigma_j^2(u)f_{j}(u)$$
for all $u\in[0,1]$ converges pointwise [and hence, uniformly] to a continuous function $\sigma_\infty^2$, bounded away from 0, as $N\to\infty$. 

\item[$\tilde{A}_4$.] With $\varepsilon_i :=Y_i-\mu(X_i)$ for all $i=1,\dots,N$, there exists $\sigma>0$ such that $\EE[|\varepsilon_i|^p|X_i=t]\leq \sigma^p$ for all $i,t$ and some $p>2$. 

\item[$\tilde{A}_5$.] The function $\mu$ is decreasing and has a continuous first derivative on $[0,1]$ such that $\inf_{u\in[0,1]}|\mu'(u)|>0$
\end{itemize}
For notational convenience, we do not make it explicit in the notation that $F_X,f_X,\sigma_X,m$ may depend on $N$.
\newline
\newline
{\bf Remark:} The pointwise convergence of $f_X$ to $f_{\infty}$ implies uniform convergence because by assumption $\tilde{A}_1$, the class of functions $\{f_j\}$ is uniformly equicontinuous, which then implies that the class $\{f_X\}$ is also uniformly equicontinuous. Also, the pointwise convergence of $\sigma_X^2$ to $\sigma_{\infty}^2$ guarantees uniform convergence, because the class of functions $\{\sigma_X^2\}$ is uniformly equicontinuous: this follows from the uniform boundedness of the class $\{f_j\}$ assumed in $\tilde{A}_0$, the uniform boundedness of $\{\sigma_j^2\}$, which is a consequence of 
$\tilde{A}_4$, and the uniform equicontinuity of the classes $\{f_j\}$ and $\{\sigma_j^2\}$ assumed in $\tilde{A}_1$.
%\newline
%\newline
%{\bf Remark for Cecile:} We want to deduce that $\sigma^2_X$ converges \emph{uniformly} to $\sigma^2_{\infty}$, provided it converges pointwise. As I argued in the e-mail, it suffices to show that the sequence $\{\sigma^2_X\}$ is uniformly equicontinuous, and this follows if in addition to $\tilde{A}_1$ we further know that $\{\sigma_j^2\}$ and $\{f_j\}$ are uniformly bounded. That $\sigma_j^2$ is uniformly bounded we can deduce from $\tilde{A}_4$, as you do in fact, later in the proof. But it is not clear to me that the $f_j$'s are uniformly bounded. Note that the $f_j$'s \emph{do not converge pointwise}, it is $f_X$ a weighted average that does. So, yes, $f_X$ converges pointwise and is uniformly equicontinuous and therefore uniformly bounded. But can we deduce the uniform boundedness of the sequence $\{f_j\}$ from this fact? If so, then we would not need $\tilde{A}_0$. 
%\newline
%\newline
\begin{theorem}\label{theo: loiU}
With $t\in(0,1)$ fixed, and $a=\mu(t)+N^{-1/3}x$ for some fixed $x\in\RR$, under Assumptions $\tilde{A}_1$ through $\tilde{A}_4$ and A4, we have
\begin{eqnarray*}
N^{1/3}(U_N(a)-g(a))\to_d \left(\frac{2\sigma_\infty(t)}{|\mu'(t)|f_\infty(t)}\right)^{2/3} \mathbb{Z}\mbox{ as } N\to\infty,
\end{eqnarray*}
where $\mathbb{Z} := \argmax_{u \in \RR} \{W(u) - u^2\}$, $W$ being a standard two-sided Brownian motion starting at 0,  has the so-called Chernoff's distribution. 
\end{theorem}
An interesting feature of the estimator $U_N$ is that its asymptotic behavior does not depend on the way the $N$ data are allocated on the different servers. The direct estimator $\hat\mu_N$ shares this feature, as is shown in the next result. 

\begin{theorem}
\label{theo:forasy} 
Under the same assumptions as in Theorem \ref{theo: loiU}, with $t\in(0,1)$ fixed, we have
\begin{eqnarray*}
N^{1/3}(\widehat\mu_N(t)-\mu(t))\to_d \left(\frac{4\sigma_\infty^2(t)|\mu'(t)|}{f_\infty^2(t)}\right)^{1/3}\mathbb{Z}\mbox{ as } N\to\infty,
\end{eqnarray*}
where $\mathbb{Z}$ is as defined in Theorem \ref{theo: loiU}. 
\end{theorem}
\noindent
{\bf Remark:} The estimators $\widehat\mu_N(t)$ and $U_N(a)$ have the \emph{same asymptotic distributions} (when centered around their respective estimands and scaled by the factor $N^{1/3}$) as the corresponding global isotonic estimators, $\hat{\mu}_{N,G}$ and $U_{N,G}$ defined in \eqref{eq:switchglobal} and \eqref{global-inv-est} respectively. In other words, the asymptotic distributions of the estimators $N^{1/3}(U_{N,G} - g(a))$ and 
$N^{1/3}(\hat{\mu}_{N,G}(t) - \mu(t))$ are those arising in Theorems \ref{theo: loiU} and \ref{theo:forasy} respectively. The limit distributions of the global estimators can be established by the same set of techniques as used in the proofs of Theorems \ref{theo: loiU} and \ref{theo:forasy}. Thus, the new estimators proposed in this paper not only circumvent the super-efficiency phenomenon but recover the asymptotic properties of their corresponding global versions. 
We note that the global isotonic estimators $\widehat\mu_{N,G}(t)$ and $U_{N,G}(a)$ also possess the \emph{uniformly bounded maximal MSE} property for their respective estimands, i.e. exact analogues of the results in Theorems \ref{theo: supereffdirect} and \ref{theo: supereffinv} hold for $N^{1/3}(U_{N,G} - g(a))$ and $N^{1/3}(\hat{\mu}_{N,G}(t) - \mu(t))$ respectively, and can be established by similar techniques as used in the proofs of these two theorems. 
\newline
\newline
{\bf Remark:} The setting of the theorems in this section with a growing sequence of sub-populations such that conditions $\tilde{A}_1$ through $\tilde{A}_5$ 
hold is not difficult to satisfy. Consider, for example, $m = \lfloor N^{1/4} \rfloor$ and $P_j$ has density $f_j(u) = (1-\epsilon_j) f_0(u) 
+ \epsilon_j f_1(u)$ where $f_0$ and $f_1$ are Lipschitz continuous densities bounded away from 0 and $\infty$ on $[0,1]$, $0 < \epsilon_j < 1$ for all $j$, the sequence $\{\epsilon_j\}$ is decreasing to 0 and $\sum_{j=1}^m \epsilon_j = o(m)$, which is easy to arrange. Let the distribution of the $X_i$'s be $P_1$ for $i = 1, 2, \ldots, \lfloor N/m \rfloor$, $P_2$ for $\lfloor N/m \rfloor + 1 \leq i \leq 2 \lfloor N/m \rfloor$, $\ldots$, and $P_m$ for 
$(m-1) \lfloor N/m \rfloor  \leq i \leq N$. For each $i$, the regression model is $Y = \mu(X_i) + \epsilon_i$ where the $\epsilon_i$'s are i.i.d. $N(0,\sigma^2)$ (say) and independent of the $X_i$'s, which are also mutually independent, and $\mu$ satisfies all the desired conditions in this manuscript, in particular 
$\tilde{A}_5$. Then, it is easy to check that all the five conditions at the beginning of this section hold, with $f_{\infty} = f_0$ and  
$\sigma_{\infty}^2(u) = \sigma^2\,f_{\infty}(u)$. 
\newline
\newline
The proof of Theorem 3.6 is in the appendix. The proof of Theorem 3.7 follows.
\newline
\newline
{\bf Proof of Theorem 3.7.} It follows from the switch relation \eqref{eq: switch} that for all fixed $t\in(0,1)$, with $a=\mu(t)+N^{-1/3}x$ we have
\begin{eqnarray*}
\PP\left(N^{1/3}(\widehat\mu_N(t)-\mu(t))<x\right)
&=&\PP\left(\widehat\mu_N(t)<\mu(t)+N^{-1/3}x\right)\\
&=&\PP\left(t>U_N(a)\right)\\
&=&\PP\left(N^{1/3}(U_N(a)-g(a))<N^{1/3}(t-g(a))\right).
\end{eqnarray*}
Now, $N^{1/3}(t-g(a))=xg'(\mu(t))+o(1)=x|\mu'(t)|^{-1}+o(1),$ so it follows from Theorem \ref{theo: loiU} that
\begin{eqnarray*}
\lim_{N\to\infty}\PP\left(N^{1/3}(\widehat\mu_N(t)-\mu(t))<x\right)=
\PP\left(\left(\frac{2\sigma_\infty(t)}{|\mu'(t)|f_X(t)}\right)^{2/3} \mathbb{Z}<\frac{x}{|\mu'(t)|}\right),
\end{eqnarray*}
using that the Chernoff distribution $\mathbb{Z}$ is continuous (see e.g. \cite{GW01}).  
\hfill{$\Box$}

\section{Discussion} 
We have proposed new estimators for distributed computing in the isotonic regression problem whose computations are not anymore onerous than that of the respective global isotonic estimators, which replicate the properties of the global estimators, and do not suffer from the superefficiency phenomenon unlike the BDSE. The key change from the BDS procedure lies in smoothing the data on local servers followed by isotonization on the central server, an `SI' (smoothing-isotonization) procedure. We note here that such `SI' procedures and their converse ('IS') procedures have been studied in monotone function problems, though \emph{not in distributed computing environments} and \emph{not under the heterogeneity setting of  our paper} . See, for example, \cite{mukerjee1988monotone},  
%does IS. Smoothness assumptions are stronger than once continuous differentiability. Normal distributions obtained under a variety of bandwidths $b_n$ with $n^{-1/3} = o(b_n)$. 
%\newline
\cite{Mammen91}, \cite{van2003smooth}, \cite{anevski2006general} and \cite{groeneboom2010maximum}. 
% does both IS and SI, asymptotic equivalence, asymptotics under two derivatives on $m$, limiting normal distributions with an $n^{2/5}$ convergence rate as in standard kernel density estimation using an $n^{-1/5}$ order bandwidth. Followed up by 
%and \cite{Groeneboom et. al. (2010)}.  
%\newline
%studied estimation of a decreasing density (under a once differentiable assumption) via an isotonized kernel estimator (SI) using a bandwidth of order $n^{-1/3}$. In this case, the limit distribution is non-Gaussian and depends on the choice of the kernel  (Theorem 2.2) and a similar phenomenon was observed in Anevski and Hossjer (2006) (Theorems 5 and 6.1). 
\newline
\newline
The ideas in this paper also have certain connections to other work in the monotone function literature which are worth mentioning. \cite{zhang2001asymptotic}  study isotonic estimation of a decreasing density with histogram-type data  based on i.i.d. data under a once differentiable assumption on the density. The domain of the density is split into bins, and the counts in each bin are available. When the number of bins grows at a rate faster than $n^{1/3}$, Theorem 4.6 of this paper shows that the isotonic estimate based on binned data recovers the Chernoff-type asymptotic distribution of the classical Grenander estimator. A similar phenomenon transpires in our problem. The $(C_{\ell k}, T_{\ell k})$ pair records the number of observations in the bin $I_k$ and the sum of the responses in that bin respectively, for the $\ell$'th server. Once these are transferred to the central server, we sum across $\ell$ to find the total number of observations in $I_k$ and the sum of the responses corresponding to all those observations and construct our isotonic estimator using these statistics. In our problem, $K$ grows faster than $N^{1/3}$ and we obtain a Chernoff limit for the pooled estimator.  

This naturally raises the question as to how the number of bins $K$ for the smoothing step on the local servers would influence the distribution of the estimators developed in this paper.  When $N^{1/3} = o(K)$, the grid is sufficiently dense and the corresponding bins sufficiently small, so that our isotonized regressogram estimator recovers the asymptotics of the classical, i.e. global isotonic regression estimator, but this will no longer be the case when $K \sim 
N^{1/3}$ or $K = o(N^{1/3})$. When $K \sim N^{1/3}$, the results of \cite{zhang2001asymptotic} (Theorem 3.3 and Corollary 4.4) and \cite{tang2012likelihood} (Theorem 3.7) who study monotone function estimation with covariates supported on a grid indicate that the limit distribution of the isotonized regressogram estimator at a point will neither be normal, nor will it be given by Chernoff's distribution. When $K = o(N^{1/3})$, the grid is sparse enough so that the regressogram estimates are ordered with probability increasing to one, so that the isotonized regressogram estimator agrees with the original estimator with increasing probability, and the results in \cite{zhang2001asymptotic} (Theorem 4.1) and \cite{tang2012likelihood} (Theorem 3.1) suggest an asymptotic normal distribution for our proposed estimator. We do not go into a full investigation of the details of these asymptotics in the distributed setting, however, since this is not relevant to the goal of the current work: produce a pooled estimator whose properties mimic the global estimator. 

We believe that similar estimators can be proposed for distributed convex regression. For convex regression, a BDS type estimator is expected to fail completely, since the global convex least squares estimator is itself asymptotically biased, as suggested by the extensive simulation experiments in  \cite{azadbakhsh2014computing}. However, a convexified regressogram estimator in the spirit of the one considered in this paper, ought to be able to recover the properties of the global convex LS estimator provided $K$ is selected appropriately: we conjecture that in the convex case $K$ should be taken to be $K^{-1} = o(N^{-1/5})$. This will provide a possible avenue for future research. 
\newline
\newline
{\bf Acknowledgement:} We would like to thank Professor Ya'acov Ritov for some fruitful discussions (with the first author) that inspired the construction of the estimates proposed in this paper.  The contribution of the second author has been conducted as part of the project Labex MME-DII (ANR11-LBX-0023-01).

\section{Appendix}
\begin{lemma}\label{lem: DKWinv} 
Assume that the distribution function $F_X$ taken from \eqref{def: FX} has a density function $f_X$ on $[0,1]$ that satisfies \eqref{eq: HfX} for some positive numbers $C_1,C_2$. 
Let $F_{N}$ be the empirical distribution function taken from \eqref{def: FN} and let $F_{N}^{-1}$ be the corresponding empirical quantile function. We then have
\begin{eqnarray}\label{eq: DKW}
\PP\left(\sup_{t\in[0,1]}\vert F_{N}(t)-F_X(t)\vert>x\right)&\leq&2\sum_{j=1}^m\exp(-2n_jx^2)
\end{eqnarray}
and
\begin{equation}\label{eq: DKWinv}
\PP\left(\sup_{t\in[0,1]}|F_{N}^{-1}(t)-F_X^{-1}(t)|>x\right)\leq 4\sum_{j=1}^m\exp(-2n_jC_1^2x^2)
\end{equation}
for all $N$ and $x>0$. 
%Moreover, for all $p>0$ there exists $K_p>0$ that depends on $c$ and $p$ only, such that for all $n$,
%\begin{equation}\label{eq: DKWinvE}
%\EE\left(\sup_{t\in[0,1]}|F_{n}^{-1}(t)-F^{-1}(t)|^p\right)\leq K_pn^{-p/2}.
%\end{equation}
\end{lemma}

\paragraph{Proof}
Let  $F_{Xj}$ denote the common distribution function of the $X_{i}$'s from sample $j$ and denote by $(X_{ji},Y_{ji})$, $i=1,\dots,n_j$ the observations from sample $j$. It follows from the  triangle inequality that
$$\sup_{t\in[0,1]}|F_N(t)-F_X(t)|\leq \sum_{j=1}^m \frac {n_j}N \sup_{t\in[0,1]}\left\vert\frac{1}{n_j}\sum_{i=1}^{n_j}\mathds{1}_{X_{ji}\leq t}-F_{Xj}(t)\right\vert$$
where we recall that $\sum_{j=1}^mn_j=N$. Hence, for all $x>0$ we have
\begin{eqnarray*}
\PP\left(\sup_{t\in[0,1]}\vert F_{N}(t)-F_X(t)\vert>x\right)&\leq&
\PP\left(\sup_{t\in[0,1]}\left\vert\frac{1}{n_j}\sum_{i=1}^{n_j}\mathds{1}_{X_{ji}\leq t}-F_{Xj}(t)\right\vert>x\mbox{ for some }j\in\{1,\dots,m\}\right)\\
&\leq&\sum_{j=1}^m \PP\left(\sup_{t\in[0,1]}\left\vert\frac{1}{n_j}\sum_{i=1}^{n_j}\mathds{1}_{X_{ji}\leq t}-F_{Xj}(t)\right\vert>x\right).
%2m\exp(-2nx^2)
\end{eqnarray*}
Since for all fixed $j$, the random variables $X_{ji},\ i=1,\dots,n_j$ are i.i.d. with distribution function $F_{Xj}$, it follows from Corollary 1 in \cite{massart1990tight} that
\begin{equation*}
\PP\left(\sup_{t\in[0,1]}\left\vert\frac{1}{n_j}\sum_{i=1}^{n_j}\mathds{1}_{X_{ji}\leq t}-F_{Xj}(t)\right\vert>x\right)\leq 2\exp(-2n_jx^2).
\end{equation*}
Combining the two preceding displays completes the proof of \eqref{eq: DKW}.

Now, consider \eqref{eq: DKWinv}. Since $f_X$ is supported on $[0,1]$, both $F_{N}^{-1}$ and  $F_X^{-1}$ take values in $[0,1]$ so the sup-distance between those functions is less than or equal to one. This means that the probability on the left hand side of  \eqref{eq: DKWinv} is equal to zero for all $x\geq 1$. Hence,  it suffices to prove \eqref{eq: DKWinv} for $x\in(0,1)$. As is customary, we use the notation $y_{+}=\max(y,0)$ and $y_{-}=-\min(y,0)$ for all real numbers $y$. This means that $|y|=\max(y_{-},y_{+})$. Recall the switching relation for the empirical distribution and empirical quantile functions: for arbitrary $a\in[0,1]$ and $t\in[0,1]$, we have \begin{equation}\label{eq: switchEmpir}
F_{N}(a)\geq t\Longleftrightarrow a\geq F_{N}^{-1}(t).
\end{equation} For all $x\in(0,1)$ we then have
\begin{eqnarray*}
\PP\left(\sup_{t\in[0,1]}(F_{N}^{-1}(t)-F_X^{-1}(t))_{+}>x\right)& = &\PP\left(\exists t\in[0,1]:\ F_{N}^{-1}(t)>x+F_X^{-1}(t)\right)\\
& = &\PP\left(\exists t\in[0,1]:\ t>F_{N}(x+F_X^{-1}(t))\right).
\end{eqnarray*}
Using $t=F_X(F_X^{-1}(t))$ together with the change of variable $u=x+F_X^{-1}(t)$ we obtain
\begin{eqnarray*}
\PP\left(\sup_{t\in[0,1]}(F_{N}^{-1}(t)-F_X^{-1}(t))_{+}>x\right)
& \leq &\PP\left(\exists u\geq x:\ F_X(u-x)>F_{N}(u)\right)\\
& = &\PP\left(\exists u\in(x,1):\ F_X(u-x)>F_{N}(u)\right).
\end{eqnarray*}
For the last equality, we use the fact that $F_X(u-x)\leq 1=F_{N}(u)$ for all $u\geq 1$, and $F_X(u-x)=0\leq F_{N}(u)$ for all $u\leq x$. With $C_1$ taken from \eqref{eq: HfX} we have 
 $F_X(u-x)< F_X(u)-C_1x$ for all $x\in(0,1)$ and $u\in(x,1)$. Combining this to the  previous display yields
 \begin{eqnarray}\notag \label{eq: DKWinv+}
\PP\left(\sup_{t\in[0,1]}(F_{N}^{-1}(t)-F_X^{-1}(t))_{+}>x\right)& \leq &\PP\left(\exists u\in(x,1):\ F_X(u)-F_{N}(u)>C_1x\right)\\
\notag & \leq &\PP\left(\sup_{u\in\RR}|F_X(u)-F_{N}(u)|>C_1x\right)\\
&\leq&2\sum_{j=1}^m\exp(-2n_jC_1^2x^2).
\end{eqnarray}
For the last inequality, we used \eqref{eq: DKW}. On the other hand, for all $x\in(0,1)$  we have
\begin{eqnarray*}
\PP\left(\sup_{t\in[0,1]}(F_{N}^{-1}(t)-F_X^{-1}(t))_{-}>x\right)& \leq &\PP\left(\exists t\in[0,1]:\ F_{N}^{-1}(t)< F_X^{-1}(t)-x\right)\\
& \leq &\PP\left(\exists u\in(x,1):\ F_{N}^{-1}(F_X(u))\leq  u-x\right),
\end{eqnarray*}
using the change of variable $u=F_X^{-1}(t)$. Hence, with the switching relation we obtain
\begin{eqnarray*}
\PP\left(\sup_{t\in[0,1]}(F_{N}^{-1}(t)-F_X^{-1}(t))_{-}>x\right)
&\leq&\PP\left(\exists u\in(x,1):\ F_X(u)\leq F_{N}(u-x)\right)\\
& \leq &\PP\left(\exists u\in(x,1):\ F_X(u-x)+C_1x < F_{N}(u-x)\right),
\end{eqnarray*}
using that $F_X(u-x)< F_X(u)-C_1x$
 for all $x\in(0,1)$ and $u\in(x,1)$.  Using again \eqref{eq: DKW} together with the change of variable $v=u-x$, we arrive at
\begin{eqnarray*}
\PP\left(\sup_{t\in[0,1]}(F_{N}^{-1}(t)-F_X^{-1}(t))_{-}>x\right)
& \leq &\PP\left(\sup_{v\in\RR}|F_X(v)-F_{N}(v)|>C_1x\right)\\
&\leq&2\sum_{j=1}^m\exp(-2n_jC_1^2x^2).
\end{eqnarray*}
Combining the previous display  with \eqref{eq: DKWinv+} completes the proof of \eqref{eq: DKWinv} since $|y|\leq y_{-}+y_{+}$ for all $y\in\RR$. %Then, \eqref{eq: DKWinvE} follows from  \eqref{fubini} combined to \eqref{eq: DKWinv}. 
\hfill{$\Box$}

\begin{lemma}\label{lem: DKWinv} 
Under the assumptions of Theorem \ref{theo: loiU},  for all $p>0$ there exists $K_p>0$  such that for all $N$,
\begin{equation}\label{eq: DKWinvE}
\EE\left(\sup_{t\in[0,1]}|F_{N}(t)-F_X(t)|^p\right)\leq K_pN^{-p/2}.
\end{equation}
\end{lemma}

\paragraph{Proof.}
It follows from  the Fubini theorem that
\begin{eqnarray*}
\EE\left(\sup_{t\in[0,1]}|F_{N}(t)-F_X(t)|^p\right)&=&\int_0^\infty
\PP\left(\sup_{t\in[0,1]}|F_{N}(t)-F_X(t)|^p>x\right)dx\\
&=&\int_0^\infty px^{p-1}
\PP\left(\sup_{t\in[0,1]}|F_{N}(t)-F_X(t)|>x\right)dx.
\end{eqnarray*}
Combining this with \eqref{eq: DKW} and the fact that a probability cannot be larger than one then yields
\begin{eqnarray*}
\EE\left(\sup_{t\in[0,1]}|F_{N}(t)-F_X(t)|^p\right)
&\leq &N^{-p/3}+2\sum_{j=1}^m\int_{N^{-1/3}}^\infty px^{p-1}\exp(-2n_jx^2)dx\\
&\leq &N^{-p/3}+2N\int_{N^{-1/3}}^\infty px^{p-1}\exp(-2N^{2/3}(\log N)^3x^2)dx
\end{eqnarray*}
for sufficiently large $N$, where we used \eqref{eq: lambdaj} for the last inequality. The result follows by computing the integral on the right-hand side. \hfill{$\Box$}
\newline
\newline
{\bf Proof of Lemma \ref{lem: EXLambda}.} For all $a\in\RR$ and $u\in\{\overline x_0,\dots,\overline x_K\}$ such that $|u-g(a)|\geq N^{-1/3}$, define
\begin{equation}\label{eq: eau}
e(a,u)=E^X\Lambda_N(u)-E^X\Lambda_N\left(\frac{[Kg(a)]}{K}\right)-a\left(F_N(u)-F_N\left(\frac{[Kg(a)]}{K}\right)\right).
\end{equation}
By definition of $\Lambda_N$ we have
\begin{equation}\notag
e(a,u)=\frac 1N\sum_{i=1}^N\mu(X_i)\left( \mathds{1}_{X_i\leq u}-\mathds{1}_{X_i\leq [Kg(a)]K^{-1}}\right)-a\left(F_N(u)-F_N\left(\frac{[Kg(a)]}{K}\right)\right).
\end{equation}
Now, $X_i\neq [Kg(a)]K^{-1}$ for all $i$, almost surely  since $X_i$ has a continuous distribution function, so \eqref{eq: Hmu} implies that
\begin{equation}\label{eq: dleau}\notag
\left\vert\mu(X_i)-\mu\left(\frac{[Kg(a)]}{K}\right)\right\vert\geq \left\vert X_i-\frac{[Kg(a)]}{K}\right\vert C_3,
\end{equation}
implying that
\begin{eqnarray}\label{eq: ftoe}
e(a,u)&\leq& \left(\mu\left(\frac{[Kg(a)]}{K}\right)-a\right)\left(F_N(u)-F_N\left(\frac{[Kg(a)]}{K}\right)\right)-C_3f(a,u)
\end{eqnarray}
with a decreasing function $\mu$, where
\begin{eqnarray}\label{eq: fau}
f(a,u)=\frac {1}{N}\sum_{i=1}^N\left(X_i-\frac{[Kg(a)]}{K}\right)\left(\mathds{1}_{X_i\leq u}-\mathds{1}_{X_i\leq [Kg(a)]K^{-1}}\right).
\end{eqnarray}
Using again \eqref{eq: Hmu}, we obtain that for all $a\in[\mu(1),\mu(0)],$
\begin{eqnarray}\label{eq: Deltamu}\notag
\left\vert\mu\left(\frac{[Kg(a)]}{K}\right)-a\right\vert&=&
\left\vert\mu\left(\frac{[Kg(a)]}{K}\right)-\mu\circ g(a)\right\vert\\
&\leq&C_4K^{-1}.
\end{eqnarray}
On the other hand, since $F_X$ has a bounded derivative that satisfies \eqref{eq: HfX} we have
\begin{eqnarray}\label{eq: DeltaF}
\notag
\left\vert F_X(u)-F_X\left(\frac{[Kg(a)]}{K}\right)\right\vert
&\leq&\left\vert F_X(u)-F_X(g(a))\right\vert+\left\vert F_X(g(a))-F_X\left(\frac{[Kg(a)]}{K}\right)\right\vert\\ \notag
&\leq&C_2\left(|u-g(a)|+K^{-1}\right)\\
&\leq &2C_2|u-g(a)|
\end{eqnarray}
for sufficiently large $N$,  using that $K^{-1}=o(N^{-1/3})$ whereas $|u-g(a)|\geq N^{-1/3}$ for the last inequality.
Next, since $m\leq N$, it follows from \eqref{eq: DKW} in the appendix that
\begin{eqnarray*}
\PP\left(\sup_{t\in[0,1]}\vert F_{N}(t)-F_X(t)\vert>x\right)&\leq&2N\exp\left(-2x^2\min_{1\leq j\leq m}n_j\right)
\end{eqnarray*} 
for all $x>0$.
With  \eqref{eq: lambdaj}, we obtain
\begin{eqnarray*}
\PP\left(\sup_{t\in[0,1]}\vert F_{N}(t)-F_X(t)\vert>x\right)&\leq&2N\exp\left(-2x^2N\lambda\right)
\end{eqnarray*} 
for all $x>0$.
With $\widetilde {\cal E}_N$ the event that 
\begin{equation}\label{eq: tildeEN}
\sup_{t\in[0,1]}\vert F_{N}(t)-F_X(t)\vert\leq C_2N^{-1/3}(\log N)^{-1}
\end{equation}
 we conclude from the previous display that
\begin{eqnarray}\label{eq: DKWO}
1-\PP(\widetilde {\cal E}_N)\notag
&\leq&2N\exp\left(-2C_2^2N^{1/3}\lambda(\log N)^{-2}\right)\\
&\ll&N^{-\theta},
\end{eqnarray} 
where we used \eqref{eq: lambdaj} for the last claim.
Combining \eqref{eq: ftoe}, \eqref{eq: Deltamu} and \eqref{eq: DeltaF} proves that  on $\widetilde{\cal E}_N$, we have
\begin{eqnarray}\label{eq: maje}\notag
e(a,u)&\leq&C_4K^{-1}\left(\left|F_X(u)-F_X\left(\frac{[Kg(a)]}{K}\right)\right|+2\sup_{t\in[0,1]}\vert F_{N}(t)-F_X(t)\vert\right)-C_3f(a,u)\\ 
&\leq&2C_2C_4K^{-1}(|u-g(a)|+N^{-1/3})-C_3f(a,u)
\end{eqnarray}
for all $a\in[\mu(1),\mu(0)]$. The inequality in \eqref{eq: ftoe} holds also for $a\not\in[\mu(1),\mu(0)]$, and in that case, 
\[
\frac{[Kg(a)]}{K}=g(a)=
\begin{cases} 0&\mbox{if } a>\mu(0)\\
1&\mbox{if } a<\mu(1),
\end{cases}
\]
implying that
$$
\left(\mu\left(\frac{[Kg(a)]}{K}\right)-a\right)\left(F_N(u)-F_N\left(\frac{[Kg(a)]}{K}\right)\right)\leq 0.$$
Hence, the inequality in \eqref{eq: maje} holds for all $a\in\RR$ and $u\in\{\overline x_0,\dots,\overline x_K\}$.
Using that $K^{-1}=o(N^{-1/3})$ whereas $|u-g(a)|\geq N^{-1/3}$, we conclude that on $\widetilde {\cal E}_N$,
\begin{eqnarray*}
e(a,u)
&\leq&o(u-g(a))^2-C_3f(a,u)
\end{eqnarray*}
uniformly over all $a$ and $u$ such that $|u-g(a)|\geq N^{-1/3}$.
 Hence, it suffices to prove that with $f(a,u)$ taken from \eqref{eq: fau}, there exists $\tilde{c}>0$ that only depends on $C_1$ such that on an event ${\cal E}_N$ whose probability is larger than $1-N^{-\theta}$, and such that ${\cal E}_N\subset\widetilde{\cal E}_N$, we have 
\begin{equation}\label{eq: pfau}
f(a,u)\geq \tilde{c}(u-g(a))^2\ for\  all\ a\in\RR,\ u\in\{\overline x_0,\dots,\overline x_K\}\ such\ that\ |u-g(a)|\geq N^{-1/3}.
\end{equation}
Similar to \eqref{eq: DKWO}, if follows from \eqref{eq: DKWinv} in the appendix that 
\begin{eqnarray}\label{eq: DKWinvO}\notag
\PP\left(\sup_{x\in[0,1]}\vert F_{N}^{-1}(x)-F_X^{-1}(x)\vert>\frac{N^{-1/3}}{\log N}\right)
&\leq&4N\exp\left(-2C_1^2N^{1/3}(\log N)^{-2}\lambda\right)\\
&\ll&N^{-\theta}.
\end{eqnarray} 
In the sequel, we consider
\begin{equation}\notag
{\cal E}_N=\widetilde{\cal E}_N\cap\left\{\sup_{x\in[0,1]}\vert F_{N}^{-1}(x)-F_X^{-1}(x)\vert\leq N^{-1/3}(\log N)^{-1}\right\}.
\end{equation}
It follows from \eqref{eq: DKWO} and \eqref{eq: DKWinvO} that $1-\PP({\cal E}_N)\ll N^{-\theta}$ so in particular, $\PP({\cal E}_N)\geq 1-N^{-\theta}$ for sufficiently large $N$. It remains to show that \eqref{eq: pfau} holds on ${\cal E}_N$. Since $X_1,\dots,X_N$ are independent with a continuous distribution function, they are all distinct from each other and for all $i$, there exists a (unique) random $j$ such that $X_i=F_N^{-1}(j/N)$, where $F_N^{-1}$ is the empirical quantile function corresponding to $X_1,\dots,X_N$. Hence, reordering the terms in the sum in \eqref{eq: fau}, we  obtain that
\begin{eqnarray*}
f(a,u)&=&\frac {1}{N}\sum_{i=1}^N\left(F_N^{-1}(i/N)-\frac{[Kg(a)]}{K}\right)\left(\mathds{1}_{F_N^{-1}(iN^{-1})\leq u}-\mathds{1}_{F_N^{-1}(iN^{-1})\leq [Kg(a)]K^{-1}}\right)\\
&=&\frac {1}{N}\sum_{i=1}^N\left(F_N^{-1}(i/N)-\frac{[Kg(a)]}{K}\right)\left(\mathds{1}_{iN^{-1}\leq F_N(u)}-\mathds{1}_{iN^{-1}\leq F_N([Kg(a)]K^{-1})} \right).
\end{eqnarray*}
Using that $F_N^{-1}$ is constant on all intervals $((i-1)N^{-1},iN^{-1}]$ we arrive at
\begin{eqnarray*}
f(a,u)&=&\int^{F_N(u)}_{F_N([Kg(a)]K^{-1})}\left(F_N^{-1}(x)-\frac{[Kg(a)]}{K}\right)dx.
\end{eqnarray*}
Hence, on ${\cal E}_N$ we have
\begin{eqnarray*}
&&\left|f(a,u)-\int^{F_N(u)}_{F_N([Kg(a)]K^{-1})}\left(F_X^{-1}(x)-\frac{[Kg(a)]}{K}\right)dx\right|\\
&&\qquad\leq\left|F_N(u)-F_N([Kg(a)]K^{-1})\right|\times \sup_{x\in[0,1]}\vert F_{N}^{-1}(x)-F_X^{-1}(x)\vert\\
%&&\qquad\leq C_2\left(|u-g(a)|+K^{-1}+2N^{-1/3}(\log N)^{-1}\right)\times \sup_{x\in[0,1]}\vert F_{N}^{-1}(x)-F_X^{-1}(x)\vert\\
&&\qquad\leq C_2\left(|u-g(a)|+K^{-1}+2N^{-1/3}(\log N)^{-1}\right)N^{-1/3}(\log N)^{-1}
\end{eqnarray*}
for all $a,u$. 
%Hence, on ${\cal E}_N$ we have
%\begin{eqnarray*}
%f(a,u)
%&=&\int^{F_N(u)}_{F_N([Kg(a)]K^{-1})}\left(F_X^{-1}(x)-\frac{[Kg(a)]}{K}\right)dx+O\left((|u-g(a)|+N^{-1/3})N^{-1/3}(\log N)^{-1}\right),
%\end{eqnarray*}
%uniformly over $a,u$. 
%On your request, I have detailed the calculations below.
Hence,
\begin{eqnarray*}
&&\left|f(a,u)
-\int^{F_X(u)}_{F_X(g(a))}\left(F_X^{-1}(x)-g(a)\right)dx\right|\\
&&\qquad \leq C_2\left(|u-g(a)|+K^{-1}+2N^{-1/3}(\log N)^{-1}\right)N^{-1/3}(\log N)^{-1}\\
&&\qquad \qquad
+\left|\int^{F_N(u)}_{F_N([Kg(a)]K^{-1})}\left(F_X^{-1}(x)-\frac{[Kg(a)]}{K}\right)dx- \int^{F_X(u)}_{F_X(g(a))}\left(F_X^{-1}(x)-g(a)\right)dx\right|
\end{eqnarray*}
It follows that
\begin{eqnarray*}
&&\left|f(a,u)
-\int^{F_X(u)}_{F_X(g(a))}\left(F_X^{-1}(x)-g(a)\right)dx\right|\\
&&\qquad \leq C_2\left(|u-g(a)|+K^{-1}+2N^{-1/3}(\log N)^{-1}\right)N^{-1/3}(\log N)^{-1}
+K^{-1}\left|{F_X(u)}-{F_X(g(a))}\right|\\
&&\qquad \qquad
+\left|\int^{F_N(u)}_{F_N([Kg(a)]K^{-1})}\left(F_X^{-1}(x)-\frac{[Kg(a)]}{K}\right)dx- \int^{F_X(u)}_{F_X(g(a))}\left(F_X^{-1}(x)-\frac{[Kg(a)]}{K}\right)dx\right|.
\end{eqnarray*}
Now,  on ${\cal E}_N$ we also have
\begin{eqnarray*}
\left|F_X^{-1}(x)-\frac{[Kg(a)]}{K}\right|&=&\left|F_X^{-1}(x)-F_X^{-1}\circ F_X\left(\frac{[Kg(a)]}{K}\right)\right|\\
&\leq&\frac 1{C_1}\left|x-F_X\left(\frac{[Kg(a)]}{K}\right)\right|\\
&\leq&\frac 1{C_1}\left(\left|F_N(u)-F_N\left(\frac{[Kg(a)]}{K}\right)\right|+C_2N^{-1/3}(\log N)^{-1}\right),
\end{eqnarray*}
for all $x$ lying between $F_N(u)$ and $F_N([Kg(a)]K^{-1})$. For such $x$'s, we obtain on ${\cal E}_N$ that
\begin{eqnarray*}
\left|F_X^{-1}(x)-\frac{[Kg(a)]}{K}\right|&\leq&\frac 1{C_1}\left(\left|F_X(u)-F_X\left(\frac{[Kg(a)]}{K}\right)\right|+3C_2N^{-1/3}(\log N)^{-1}\right)\\
&\leq&\frac {3C_2}{C_1}\left(|u-g(a)|+K^{-1}+N^{-1/3}(\log N)^{-1}\right)
\end{eqnarray*}
for all $a$ and $u$, for sufficiently large $N$. Therefore, with $K\geq 1$ we obtain on ${\cal E}_N$ that
\begin{eqnarray*}
&&\left|f(a,u)
-\int^{F_X(u)}_{F_X(g(a))}\left(F_X^{-1}(x)-g(a)\right)dx\right|\\
&&\qquad \leq 2C_2\left(|u-g(a)|+K^{-1}+2N^{-1/3}(\log N)^{-1}\right)N^{-1/3}(\log N)^{-1}
\\ 
&&\qquad\qquad +\frac {3C_2}{C_1}\left(|u-g(a)|+K^{-1}+N^{-1/3}(\log N)^{-1}\right)\left(2\sup_{u\in[0,1]}|F_N(u)-F_X(u)|+C_2K^{-1})\right)\\
&&\qquad=O\left(|u-g(a)|+K^{-1}+N^{-1/3}(\log N)^{-1}\right)\left(N^{-1/3}(\log N)^{-1}+K^{-1}\right)
\end{eqnarray*} 
on ${\cal E}_N$,  uniformly over $a\in\RR$ and $u\in\{\overline x_0,\dots,\overline x_K\}$. Now, 
 we can do the change of variable $t=F_X^{-1}(x)$ to get on ${\cal E}_N$ that
\begin{eqnarray}\label{eq: equivf}\notag
f(a,u)
&=&\int^{u}_{g(a)}\left(t-g(a)\right)f_X(t)dt\\
&&\qquad+O\left(|u-g(a)|+K^{-1}+N^{-1/3}(\log N)^{-1}\right)\left(N^{-1/3}(\log N)^{-1}+K^{-1}\right)
\end{eqnarray}
uniformly over $a\in\RR$ and $u\in\{\overline x_0,\dots,\overline x_K\}$. Here,
\begin{eqnarray}\notag
\int^{u}_{g(a)}\left(t-g(a)\right)f_X(t)dt
&\geq &C_1\int^{u}_{g(a)}\left(t-g(a)\right)dt
\end{eqnarray}
where $C_1$ is taken from \eqref{eq: HfX}, for all $a,u$. Since it is assumed that $K^{-1}=o(N^{-1/3})$, we conclude that on ${\cal E}_N$,
\begin{eqnarray*}
f(a,u)
&\geq &\frac {C_1}2(u-g(a))^2+o((g(a)-u)^2),
\end{eqnarray*}
where the small $o$-term is uniform over all $u$ and $a$ such that $|u-g(a)|\geq N^{-1/3}$. 
\newline
Hence, \eqref{eq: pfau} holds on ${\cal E}_N$  provided that $
\tilde{c}< C_1/2$ and $N$ is sufficiently large. It follows that on $\mathcal{E}_N$, for all sufficiently large $N$, 
\[e(a,u) \leq  o((u-g(a))^2 - C_3 \tilde{c}(g(a)-u)^2 \,\] 
where in view of the above proof, the small-o term can be chosen of the form 
$$ o((u-g(a))^2=2C_2C_4K^{-1}(|u-g(a)|+N^{-1/3}).$$
Therefore, for any $c < C_3 \tilde{c}$, for all sufficiently large $N$, $e(a,u) \leq -c (g(a)-u)^2$ on $\mathcal{E}_N$. This completes the proof of the lemma. \hfill{$\Box$}
\newline
\newline
{\bf Proof of Lemma \ref{lem: tail boundVN}.} For all $a\not\in[\mu(1),\mu(0)]$ and $u\in\{\overline x_0,\dots,\overline x_K\}$, define $e(a,u)$ as in \eqref{eq: eau}. We then have \eqref{eq: ftoe} where 
\[
\frac{[Kg(a)]}{K}=g(a)=
\begin{cases} 0&\mbox{if } a>\mu(0)\\
1&\mbox{if } a<\mu(1).
\end{cases}
\]
and $f$ is given by \eqref{eq: fau}. Note that \eqref{eq: Deltamu} is no longer true for $a\not\in[\mu(1),\mu(0)]$ since in such a case, $a\neq \mu\circ g(a).$ Instead, we will use
\[
\left(\mu\left(\frac{[Kg(a)]}{K}\right)-a\right)\left( F_N(u)-F_N\left(\frac{[Kg(a)]}{K}\right)\right)=
\begin{cases} (\mu(0)-a)F_N(u)&\mbox{if } a>\mu(0)\\
(\mu(1)-a)(F_N(u)-1)&\mbox{if } a<\mu(1),
\end{cases}
\]
using that $F_N(0)=0$ and $F_N(1)=1$. Since $f(a,u)\geq 0$ for all $a,u$ \eqref{eq: ftoe} yields
\begin{equation}\label{eq: majeHS}
e(a,u)\leq
\begin{cases} (\mu(0)-a)F_N(u)&\mbox{if } a>\mu(0)\\
(\mu(1)-a)(F_N(u)-1)&\mbox{if } a<\mu(1)
\end{cases}
\end{equation}
for all $a\not\in[\mu(1),\mu(0)]$ and $u\in\{\overline x_0,\dots,\overline x_K\}$.

Since $\Lambda_N(\overline x_0)-a\widetilde F_N(\overline x_0)=0$, it follows from the definition of $V_N$ that the following inequalities hold for all $x>0$ and $a>\mu(0)$:
\begin{eqnarray*}
\PP\left(V_N(a)\geq x\right)
&\leq&\PP\left(\max_{u\in\{ \widetilde F_N(\overline x_0),\dots, \widetilde F_N(\overline x_K)\},\ u\geq x}\{\Lambda_N\circ\widetilde F_N^{-1}(u)-au\}\geq 0\right)\\
&=&\PP\left(\max_{u\in\{ \widetilde F_N(\overline x_0),\dots, \widetilde F_N(\overline x_K)\},\ u\geq x}\{M_N\circ\widetilde F_N^{-1}(u)+e(a,\widetilde F_N^{-1}(u))\}\geq 0\right)
\end{eqnarray*}
where $M_N(u)=\Lambda_N(u)-\EE^X(\Lambda_N(u))$ takes the form \eqref{eq: defMN}. The first inequality in \eqref{eq: majeHS} then yields
\begin{eqnarray*}
\PP\left(V_N(a)\geq x\right)
&\leq&\PP\left(\max_{u\in\{ \widetilde F_N(\overline x_0),\dots, \widetilde F_N(\overline x_K)\},\ u\geq x}\{M_N\circ\widetilde F_N^{-1}(u)+(\mu(0)-a)u\}\geq 0\right)\\
&\leq&\sum_{k\geq 0}\PP\left(\max_{u\in\{ \widetilde F_N(\overline x_0),\dots, \widetilde F_N(\overline x_K)\},\ u\in[x2^k,x2^{k+1}]}\{M_N\circ\widetilde F_N^{-1}(u)\}\geq (a-\mu(0))x2^k\right).
\end{eqnarray*}
Let $p\geq 2$ and $\sigma>0$ such that $\EE[\varepsilon_i^p|X_i]\leq \sigma^p$ for all $i$, almost surely.
The process $M_n$ is a centered martingale under $\PP^X$ which, according to Theorem 3 in \cite{rosenthal1970subspaces}, satisfies 
\begin{eqnarray*}\label{eq: varMN}
\EE^X\left\vert M_N(u)\right\vert^p &\leq&
\frac {A_p}{N^p}\max\left\{\sum_{i=1}^N\EE^X|\varepsilon_i|^p\mathds{1}_{X_i\leq u};
\left(\sum_{i=1}^N\EE^X|\varepsilon_i|^2\mathds{1}_{X_i\leq u}\right)^{p/2}\right\}\\
&\leq&
\frac {A_p\sigma^p}{N^{p}}\max\left\{NF_N(u);
\left(NF_N(u)\right)^{p/2}\right\}\\
&\leq&
\frac {A_p\sigma^pF_N(u)}{N^{p/2}}
\end{eqnarray*}
for all $u\in[0,1]$ and $A_p=(p/2)^{p/2}2^{p+p^2/4}$. For the penultimate inequality, we used that $\EE^X|\varepsilon_i|^2\leq (\EE^X|\varepsilon_i|^p)^{2/p}$ thanks to the Holder inequality whereas for the last inequality, we used that $N\leq N^{p/2}$ and $F_N^{p/2}(u)\leq F_N(u)$. Combining the two preceding displays with  the Doob inequality yields  that for all $x>0$,
\begin{eqnarray*}
\PP\left(V_N(a)\geq x\right)
&\leq&\sum_{k\geq 0}\EE\left[\PP^X\left(\max_{u\in\{ \widetilde F_N(\overline x_0),\dots, \widetilde F_N(\overline x_K)\},\ u\in[x2^k,x2^{k+1}]}\{M_N\circ\widetilde F_N^{-1}(u)\}\geq (a-\mu(0))x2^k\right)\right]\\
&\leq&\sum_{k\geq 0}\EE\left[ \frac{A_p\sigma^pF_N(x2^{k+1})}{N^{p/2}(a-\mu(0))^p(x2^k)^p}\right].
\end{eqnarray*}
With $C_2$ taken from \eqref{eq: HfX} we conclude that for all $x>0$,
\begin{eqnarray*}
\PP\left(V_N(a)\geq x\right)
&\leq&\sum_{k\geq 0}\frac{A_p\sigma^pF_X(x2^{k+1})}{N^{p/2}(a-\mu(0))^p(x2^k)^p}\\
&\leq&\sum_{k\geq 0}\frac{2A_pC_2\sigma^p}{N^{p/2}(a-\mu(0))^p(x2^k)^{p-1}}.
\end{eqnarray*}
Since $C:=2A_pC_2\sigma^p\sum_{k\geq 0}2^{-k(p-1)}$ is finite, we conclude that 
\begin{eqnarray*}
\PP\left(V_N(a)\geq x\right)\leq\frac{C}{N^{p/2}(a-\mu(0))^px^{p-1}},
\end{eqnarray*}
which proves the first assertion. For the second assertion, since $\overline x_K=\widetilde F_N(\overline x_K)=1$, we write for $a<\mu(1)$ and $x>0$:
\begin{eqnarray*}
&&\PP\left(1-V_N(a)\geq x\right)\\
&&\qquad\leq\PP\left(\max_{u\in\{ \widetilde F_N(\overline x_0),\dots, \widetilde F_N(\overline x_K)\},\ 1-u\geq x}\{\Lambda_N\circ\widetilde F_N^{-1}(u)-au\}\geq \Lambda_N(1)-a\right)\\
&&\qquad=\PP\left(\max_{u\in\{ \widetilde F_N(\overline x_0),\dots, \widetilde F_N(\overline x_K)\},\ 1-u\geq x}\{M_N\circ\widetilde F_N^{-1}(u)-M_N(1)+e(a,\widetilde F_N^{-1}(u))\}\geq 0\right).
\end{eqnarray*}
The first inequality in \eqref{eq: majeHS} then yields
\begin{eqnarray*}
&&\PP\left(1-V_N(a)\geq x\right)\\
&&\qquad\leq\PP\left(\max_{u\in\{ \widetilde F_N(\overline x_0),\dots, \widetilde F_N(\overline x_K)\},\ 1-u\geq x}\{M_N\circ\widetilde F_N^{-1}(u)-M_N(1)-(\mu(1)-a)(1-u)\}\geq 0\right)\\
&&\qquad\leq\sum_{k\geq 0}\PP\left(\max_{u\in\{ \widetilde F_N(\overline x_0),\dots, \widetilde F_N(\overline x_K)\},\ 1-u\leq x2^{k+1}}\{M_N\circ\widetilde F_N^{-1}(u)-M_N(1)\}\geq (\mu(1)-a)x2^k\right),
\end{eqnarray*}
and we use the Doob inequality, similar as above. Details are omitted.
\hfill{$\Box$}
\newline
\newline
{\bf Proof of Theorem \ref{theo: loiU}:}  It follows from \eqref{def: UN} together with Lemma \ref{lem: tail boundUN} that with probablity tending to one,
\begin{eqnarray*}
N^{1/3}(U_N(a)-g(a))=\argmax_{u\in H_N}\{\Lambda_N(g(a)+N^{-1/3}u)-aF_N(g(a)+N^{-1/3}u)\}
\end{eqnarray*}
where $H_N$ is the set of all $u\in\RR$ such that $g(a)+N^{-1/3}u\in\{\overline x_0,\dots,\overline x_K\}$ and $|u|\leq v_N$, where $v_N$ is an arbitrary sequence that diverges to infinity as $N\to\infty$. In the sequel, we consider a sequence $v_N$ such that $v_N\leq \log N$ for all $N$. Hence, with probablity that tends to one we have
\begin{eqnarray}
&&N^{1/3}(U_N(a)-g(a))\\ \notag
&&\qquad =\argmax_{u\in H_N}\left\{N^{2/3}\left(M_N(g(a)+N^{-1/3}u)-M_N\left(\frac{[Kg(a)]}{K}\right)\right)+N^{2/3}e(a,g(a)+N^{-1/3}u)\right\}
\end{eqnarray}
where $M_N(u)=\Lambda_N(u)-E^X(\Lambda_N(u))$ for all $u\in\{\overline x_0,\dots,\overline x_K\}$ and $e$ is taken from \eqref{eq: eau}, that is
\begin{eqnarray}\notag\label{eq: elim}
e(a,g(a)+N^{-1/3}u)&=&\frac 1N\sum_{i=1}^N\mu(X_i)\left( \mathds{1}_{X_i\leq g(a)+N^{-1/3}u}-\mathds{1}_{X_i\leq [Kg(a)]K^{-1}}\right)\\
&&\qquad-a\left(F_N(g(a)+N^{-1/3}u)
-F_N\left(\frac{[Kg(a)]}{K}\right)\right).
\end{eqnarray}
We extend $M_N$ and $e(a,\ .\ )$ as constant functions in between two consecutive points in $H_N$ so that
\begin{eqnarray}\label{eq: N1/3UN}
&&N^{1/3}(U_N(a)-g(a))\notag\\ \notag
&&\qquad =\argmax_{|u|\leq v_N}\left\{N^{2/3}\left(M_N(g(a)+N^{-1/3}u)-M_N\left(\frac{[Kg(a)]}{K}\right)\right)+N^{2/3}e(a,g(a)+N^{-1/3}u)\right\}\\ 
&&\qquad\qquad+o_p(1).
\end{eqnarray}
Now, since $a=\mu(t)+N^{-1/3}x$ for some fixed $x\in\RR$ and $t\in(0,1)$, and $g'=1/\mu'\circ g$ on $(\mu(1),\mu(0))$ is bounded by assumption, we have
\begin{equation}\label{eq: g(a)-t}
g(a)=t+O(N^{-1/3}).
\end{equation}
Hence, for sufficiently large $N$, every  $X_i$ that lies between $ [Kg(a)]K^{-1}$ and $g(a)+N^{-1/3}u$ for some $|u|\leq v_N$ also lies in
$[t-N^{-1/3}\log N,t+N^{-1/3}\log N]$. This implies that for all such   $X_i$'s  there exists $\theta_i\in[t-N^{-1/3}\log N,t+N^{-1/3}\log N]$  such that
\begin{eqnarray}\label{eq: dleau}\notag
\mu(X_i)&=&\mu\left(\frac{[Kg(a)]}{K}\right)+\left(X_i-\frac{[Kg(a)]}{K}\right)\mu'(\theta_i)\\
&=&\mu\left(\frac{[Kg(a)]}{K}\right)+\left(X_i-\frac{[Kg(a)]}{K}\right)\left(\mu'(t)+o(1)\right)
\end{eqnarray}
where the small $o$-term  is uniform, by continuity of $\mu'$ over the compact interval $[t-N^{-1/3}\log N,t+N^{-1/3}\log N]$. Plugging this in \eqref{eq: elim}, and using the notation $f$ in \eqref{eq: fau}, yields
\begin{eqnarray}\notag
e(a,g(a)+N^{-1/3}u)&=&\left(\mu'(t)+o(1)\right)f(a,g(a)+N^{-1/3}u)\\
&&\quad+\left(\mu\left(\frac{[Kg(a)]}{K}\right)-a\right)\left(F_N(g(a)+N^{-1/3}u)
-F_N\left(\frac{[Kg(a)]}{K}\right)\right).
\end{eqnarray}
It can be seen from  the proof of Lemma \ref{lem: EXLambda} that \eqref{eq: tildeEN} holds on the event ${\cal E}_N$, whose probability tends to one as $N\to\infty$, implying that
\begin{eqnarray*}
F_N(g(a)+N^{-1/3}u)
-F_N\left(\frac{[Kg(a)]}{K}\right)&=&
F_X(g(a)+N^{-1/3}u)
-F_X\left(\frac{[Kg(a)]}{K}\right)+O_p(N^{-1/3}(\log N)^{-1})\\
&=&O_p(N^{-1/3}v_N+K^{-1}+N^{-1/3}(\log N)^{-1})\\
&=&O_p(N^{-1/3}v_N)
\end{eqnarray*}
uniformy over $u\in H_N$. Since \eqref{eq: Deltamu} holds for all $a\in[\mu(1),\mu(0)],$ combining the two preceding displays yields
\begin{eqnarray*}\notag
e(a,g(a)+N^{-1/3}u)&=&\left(\mu'(t)+o(1)\right)f(a,g(a)+N^{-1/3}u)+O_p(K^{-1}N^{-1/3}v_N).
\end{eqnarray*}
Next, we invoke \eqref{eq: equivf}, that holds on the event ${\cal E}_N$ uniformly over $a$ and $u$, to conclude that 
\begin{eqnarray*}\notag
e(a,g(a)+N^{-1/3}u)&=&\left(\mu'(t)+o(1)\right)\int^{g(a)+N^{-1/3}u}_{g(a)}\left(z-g(a)\right)f_X(z)dz+o_p(N^{-2/3})
\end{eqnarray*}
uniformly over $u\in H_N$, provided that $v_N\ll \min\{\log N; N^{-1/3} K\}$. By assumption, $N^{-1/3} K$ diverges to infinity as $N\to\infty$, so we can find a sequence $v_N$ that satisfies the above condition and that diverges to infinity as $N\to\infty$, as required in the definition of $H_N$. In the sequel, we consider a sequence $v_N$ that satisfies the above conditions and in addition, the below condition: 
\begin{eqnarray*}\label{eq: condvN}
v_N\ll \left(\max\left\{\sup_{|z-t|\leq N^{-1/3}\log N}|f_X(z)-f_\infty(z)|,\sup_{|z-t|\leq N^{-1/3}\log N}|f_\infty(t)-f_\infty(z)|\right\}\right)^{-1/2}.
\end{eqnarray*}
Note that by assumption, the right-hand side of the inequality in the  above display diverges to infinity as $N\to\infty$, which ensures existence of such a sequence $v_N$. We then have 
\begin{eqnarray*}\notag
e(a,g(a)+N^{-1/3}u)&=&\left(\mu'(t)+o(1)\right)\int^{g(a)+N^{-1/3}u}_{g(a)}\left(z-g(a)\right)f_\infty(z)dz+o_p(N^{-2/3}),
\end{eqnarray*}
using that for $u\geq 0$ (and similarly for $u\leq 0$),
\begin{eqnarray*}\notag
\left|\int^{g(a)+N^{-1/3}u}_{g(a)}\left(z-g(a)\right)(f_X(z)-f_\infty(z))dz\right|
&\leq& \int^{g(a)+N^{-1/3}u}_{g(a)}\left(z-g(a)\right)|f_X(z)-f_\infty(z)|dz\\
&\leq&\frac{N^{-2/3}v_N^2}{2}\sup_{|z-g(a)|\leq N^{-1/3} v_N}|f_X(z)-f_\infty(z)|
\end{eqnarray*}
uniformly for all $|u|\leq v_N$, which implies
\begin{eqnarray*}\notag
\left|\int^{g(a)+N^{-1/3}u}_{g(a)}\left(z-g(a)\right)(f_X(z)-f_\infty(z))dz\right|
&\leq&\frac{N^{-2/3}v_N^2}{2}\sup_{|z-t|\leq N^{-1/3} \log N}|f_X(z)-f_\infty(z)|\\
&=&o(N^{-2/3})
\end{eqnarray*}
thanks to \eqref{eq: g(a)-t}, \eqref{eq: condvN} and the assumption that $v_n\ll\log N$. Similarly, 
\begin{eqnarray*}\notag
\left|\int^{g(a)+N^{-1/3}u}_{g(a)}\left(z-g(a)\right)(f_\infty(z)-f_\infty(t))dz\right|
&\leq&\frac{N^{-2/3}v_N^2}{2}\sup_{|z-g(a)|\leq N^{-1/3} v_N}|f_\infty(z)-f_\infty(t)|\\
&\leq&\frac{N^{-2/3}v_N^2}{2}\sup_{|z-t|\leq N^{-1/3} \log N}|f_\infty(z)-f_\infty(t)|\\
&=&o(N^{-2/3})
\end{eqnarray*}
and therefore,
\begin{eqnarray*}\notag
e(a,g(a)+N^{-1/3}u)&=&\left(\mu'(t)+o(1)\right)\int^{g(a)+N^{-1/3}u}_{g(a)}\left(z-g(a)\right)f_\infty(z)dz+o_p(N^{-2/3})\\
&=&\left(\mu'(t)+o(1)\right)f_\infty(t)\int^{g(a)+N^{-1/3}u}_{g(a)}\left(z-g(a)\right)dz+o_p(N^{-2/3}).
\end{eqnarray*}
Hence we obtain
\begin{eqnarray}\label{eq: lime}
N^{2/3}e(a,g(a)+N^{-1/3}u)
&=&-(|\mu'(t)|+o(1))f_\infty(t)\frac{u^2}2+o_p(1).
\end{eqnarray}

On the other hand, with
\begin{eqnarray*}
Z_N(u)=N^{2/3}\left(M_N(g(a)+N^{-1/3}u)-M_N\left(\frac{[Kg(a)]}{K}\right)\right);
\end{eqnarray*}
where $M_N$ is as defined in \eqref{eq: defMN} for all $u\in\{\overline x_0,\dots,\overline x_K\}$, 
we have
\begin{eqnarray}\label{eq: ZN}
Z_N(u)&=&N^{-1/3}\sum_{j=1}^m\sum_{i=1}^{n_j}\varepsilon_{ji}\left(\mathds{1}_{X_{ji}\leq g(a)+N^{-1/3}u}-\mathds{1}_{X_{ji}\leq [Kg(a)]K^{-1}}\right)
\end{eqnarray}
where we denote by $(X_{ji},Y_{ji})$, $i=1,\dots,n_j$ the observations from sample $j$, for $j=1,\dots,m$, and $\varepsilon_{ji}=Y_{ji}-\mu(X_{ji})$. Note that the process $Z_N$ is centered and has been extended to $\RR$ by being constant in between two consecutive points in $H_N$. 
For all $u\geq v\geq 0$ in $H_N$ we have
\begin{eqnarray*}
\EE\left[Z_N(u)Z_N(v)\right]&=&N^{-2/3}\sum_{j=1}^m\sum_{i=1}^{n_j}\EE\left[\varepsilon_{ji}^2
\mathds{1}_{[Kg(a)]K^{-1}<X_{ji}\leq g(a)+N^{-1/3}u}\mathds{1}_{[Kg(a)]K^{-1}<X_{ji}\leq g(a)+N^{-1/3}v}\right]\\
&=&N^{-2/3}\sum_{j=1}^m\sum_{i=1}^{n_j}\EE\left[\sigma_j^2(X_{ji})
\mathds{1}_{[Kg(a)]K^{-1}<X_{ji}\leq g(a)+N^{-1/3}v}\right],
\end{eqnarray*}
where the last equality is obtained by conditioning with respect to $X_{ji}$ and using that $u\geq v\geq 0$. With $u,v$ fixed, this implies that
\begin{eqnarray*}
\EE\left[Z_N(u)Z_N(v)\right]
&=&N^{-2/3}\sum_{j=1}^m\sum_{i=1}^{n_j}\EE\left[\sigma_j^2(t)
\mathds{1}_{[Kg(a)]K^{-1}<X_{ji}\leq g(a)+N^{-1/3}v}\right]+o(1)\\
\end{eqnarray*}
using that for $u,v\in H_N$
%$\omega(\delta)\to0$ as $N\to\infty$ together with Lemma \ref{lem: DKWinv}, according to which,
%\begin{eqnarray*}
%&&N^{-2/3}\EE\left|\sum_{j=1}^m\sum_{i=1}^{n_j}\left(
%\mathds{1}_{[Kg(a)]K^{-1}<X_{ji}\leq g(a)+N^{-1/3}v}-\EE\left[
%\mathds{1}_{[Kg(a)]K^{-1}<X_{ji}\leq g(a)+N^{-1/3}v}\right]\right)\right|\\
%&&\qquad = N^{1/3}\EE\left|F_N(g(a)+N^{-1/3}v)-F_N(g(a))-(F_X(g(a)+N^{-1/3}v)-F_X(g(a)))\right|\\
%&&\qquad \leq 2N^{1/3}\EE\left(\sup_{u\in[0,1]}|F_N(u)-F_X(u)|\right)\\
%&&\qquad = o(1)
%\end{eqnarray*}
%and
\begin{eqnarray*}
&&\left|\EE\left[Z_N(u)Z_N(v)\right]-N^{-2/3}\sum_{j=1}^m\sum_{i=1}^{n_j}\EE\left[
\sigma_j^2(t)\mathds{1}_{[Kg(a)]K^{-1}<X_{ji}\leq g(a)+N^{-1/3}v}\right]\right|\\
&&\qquad \leq N^{-2/3}\sum_{j=1}^m\sum_{i=1}^{n_j}
\EE\left[|\sigma_j^2(X_{ji})-\sigma_j^2(t)|\mathds{1}_{[Kg(a)]K^{-1}<X_{ji}\leq g(a)+N^{-1/3}v}
\right]
\\
&&\qquad\leq N^{-2/3}\omega(N^{-1/3}\log N)\sum_{j=1}^m\sum_{i=1}^{n_j}\EE\left[
\mathds{1}_{[Kg(a)]K^{-1}<X_{ji}\leq g(a)+N^{-1/3}v}\right]
\end{eqnarray*}
where $\omega(\delta)\to0$ as $\delta\to0$ by assumption, and
\begin{eqnarray*}
N^{-2/3}\sum_{j=1}^m\sum_{i=1}^{n_j}\EE\left[
\mathds{1}_{[Kg(a)]K^{-1}<X_{ji}\leq g(a)+N^{-1/3}v}\right]
&=& N^{1/3}\left|F_X(g(a)+N^{-1/3}v)-F_X([Kg(a)]K^{-1})\right|\\
&=&O(1).
\end{eqnarray*}
Hence,
\begin{eqnarray*}
\EE\left[Z_N(u)Z_N(v)\right]
&=&N^{-2/3}\sum_{j=1}^m\sum_{i=1}^{n_j}\sigma_j^2(t)\PP\left(
[Kg(a)]K^{-1}<X_{ji}\leq g(a)+N^{-1/3}v\right)+o(1)\\
&=&N^{-2/3}\sum_{j=1}^m n_j\sigma_j^2(t)
\int_{[Kg(a)]K^{-1}}^{g(a)+N^{-1/3}v}f_j(z)dz+o(1)
\end{eqnarray*}
for all fixed real numbers $u\geq v\geq 0$. It follows that 
\begin{eqnarray*}
&&\left|\EE\left[Z_N(u)Z_N(v)\right]
-N^{-2/3}\sum_{j=1}^m n_j\sigma_j^2(t)
\int_{[Kg(a)]K^{-1}}^{g(a)+N^{-1/3}v}f_j(t)dz\right|\\
&&\qquad\leq N^{-2/3}\sum_{j=1}^m n_j\sigma_j^2(t)\omega(N^{-1/3}\log N)\left(N^{-1/3}v+O(K^{-1})\right)+o(1)\\
&&\qquad\leq o(1)N^{-1}\sum_{j=1}^m n_j\sigma_j^2(t)+o(1),
\end{eqnarray*}
since $\omega(\delta)\to0$ as $\delta\to0$. The Jensen inequality for conditional expectation combined with Assumption $\tilde{A}_4$ shows that $\sigma_j^2(t)\leq\sigma^2$ for all $i$ and $t$ and therefore, $N^{-1}\sum_{j=1}^m n_j\sigma_j^2(t)\leq \sigma^2$. This implies that
\begin{eqnarray*}
\EE\left[Z_N(u)Z_N(v)\right]
&=&N^{-2/3}\sum_{j=1}^m n_j\sigma_j^2(t)f_j(t)(N^{-1/3}v+o(N^{-1/3}))+o(1)\\
&=& \sigma_X^2(t)v+o(1).
\end{eqnarray*}
We conclude that for all $u\geq v\geq 0$,  $\EE\left[Z_N(u)Z_N(v)\right]=cov(Z_N(u),Z_N(v))$ converges to $\sigma_\infty^2(t)v$. The case of negative $u$ and $v$ can be treated likewise and therefore,
$cov(Z_N(u),Z_N(v))$ converges to $\sigma_\infty^2(t)(|u|\wedge |v|)$ if $uv\geq 0$. It can be seen similarly that it converges to zero if $uv<0$ (hence $u$ and $v$ have different signs). Hence, the covariance converges to $\sigma_\infty(t)cov(W(u),W(v))$, so we conclude from the Lindeberg-Feller theorem that jointly, 
\begin{eqnarray}\label{eq: fidis}
(Z_N(u_1),\dots,Z_N(u_k))\to_d\sigma_\infty(t)(W(u_1),\dots,W(u_k))
\end{eqnarray}
 for all fixed $u_1,\dots,u_k\in\RR$, as $N\to\infty$. Now, consider the restriction of $Z_N$ to the compact interval $[-M,M]$, for a fixed $M>0$. For all $\delta>0$ and $\epsilon>0$ we have
\begin{eqnarray}
\label{rosenapp} 
\PP\left(\sup_{|t-s|\leq\delta ;\ s,t\in[-M,M]}|Z_N(s)-Z_N(t)|\geq\epsilon\right)\leq \sum_{k=-M[\delta^{-1}]-1}^{M[\delta^{-1}]}\PP\left(2\sup_{|t-k\delta|\leq2\delta}|Z_N(k\delta)-Z_N(t)|\geq\epsilon\right).
\end{eqnarray} 
Let $\pi$ be the permutation such that the $X_{\pi(j)}$ are ordered in $j$, that is $X_{\pi(1)}<\dots<X_{\pi(N)}$ a.s. Let $\mathbb{P}_X$ denote the conditional probability given $X_1,\dots,X_N$. Since $\varepsilon_{\pi(1)},\dots,\varepsilon_{\pi(N)}$ are centered and independent under $\PP_X$,  the process $\{Z_N(k\delta)-Z_N(t),\ t\geq k\delta\}$ is a forward centered martingale whereas $\{Z_N(k\delta)-Z_N(t),\ t\leq k\delta\}$ is a reverse centered martingale conditionally on $X_1,\dots,X_N$, for all $k$. Hence,  it follows from the Doob inequality that for all $k$, 
\newline
\begin{eqnarray}\label{eq: doob}\notag
&&\PP\left(2\sup_{|t-k\delta|\leq2\delta}|Z_N(k\delta)-Z_N(t)|\geq\epsilon\right)\\
&&\qquad\qquad\leq \frac{2^p}{\epsilon^p}\left(\EE|Z_N(k\delta)-Z_N((k-2)\delta)|^p+\EE|Z_N(k\delta)-Z_N((k+2)\delta)|^p\right).
\end{eqnarray}
Note that the inequalities above are first obtained for the conditional probabilities and then integrated over the distribution of $X$ for the unconditional. Now, it follows from the Rosenthal inequality, see \cite{rosenthal1970subspaces}, that for all $k$ and a constant $C$ that depends only on $p$, we have
\begin{eqnarray*}
\EE|Z_N(k\delta)-Z_N((k+2)\delta)|^p\leq CN^{-p/3}\left(\sum_{i=1}^N\EE(|\varepsilon_i|^p \mathds{1}_{X_{i}\in I_k})+\left(\sum_{i=1}^N\EE(|\varepsilon_i|^2 \mathds{1}_{X_{i}\in I_k})\right)^{p/2}\right).
\end{eqnarray*}
Here, $I_k=(g(a)+N^{-1/3}k\delta,g(a)+N^{-1/3}(k+2)\delta]$ (at least if $g(a)+N^{-1/3}k\delta,g(a)$ and $N^{-1/3}(k+2)\delta$ both belong to $H_N$) and $p$ is taken from Assumption $\tilde A_4$. Hence, with $f_X$ taken from \eqref{eq: fX} we have
\begin{eqnarray*}
\EE|Z_N(k\delta)-Z_N((k+2)\delta)|^p&\leq &C\sigma^pN^{-p/3}\left(\sum_{i=1}^N\EE(\mathds{1}_{X_{i}\in I_k})+\left(\sum_{i=1}^N\EE( \mathds{1}_{X_{i}\in I_k})\right)^{p/2}\right)\\
&=&C\sigma^pN^{-p/3}\left(N\int_{I_k}f_X(u)du+N^{p/2}\left[\int_{I_k}f_X(u)du\right]^{p/2}\right).
\end{eqnarray*}
It follows from the Assumption $\tilde A_1$ that $f_X$ is bounded by a constant $A$ that does not depend on $N$ and therefore,
\begin{eqnarray*}
\EE|Z_N(k\delta)-Z_N((k+2)\delta)|^p&\leq &C\sigma^pN^{-p/3}\left(2AN^{2/3}\delta+\left[2AN^{2/3}\delta\right]^{p/2}\right)(1+o(1))\\
&\leq &2C\sigma^pN^{-p/3}\left[2AN^{2/3}\delta\right]^{p/2}\\
&=&2C\sigma^p\left[2A\delta\right]^{p/2}
\end{eqnarray*}
for $N$ sufficiently large. Arguing similarly for $\EE|Z_N(k\delta)-Z_N((k-2)\delta)|^p$ we conclude from \eqref{eq: doob} that there exists $C>0$ that depends only on $p$ and $A$ such that
\begin{eqnarray*}
\PP\left(2\sup_{|t-k\delta|\leq2\delta}|Z_N(k\delta)-Z_N(t)|\geq\epsilon\right)&\leq&  C\sigma^p \epsilon^{-p}\delta^{p/2}
\end{eqnarray*}
for all $k$. Summing up this inequality over all $k$ on the right-hand side of \eqref{rosenapp}, we obtain that there exists $C>0$ that depends only on $p$ and $A$ such that for all $\delta>0$ and $\epsilon>0$,
\begin{eqnarray*}
\PP\left(\sup_{|t-s|\leq\delta\ ;\ s,t\in[-M,M]}|Z_N(s)-Z_N(t)|\geq\epsilon\right)\leq CM\sigma^p\epsilon^{-p}\delta^{-1+p/2}.
\end{eqnarray*}
Since $p>2$, this converges to zero as $\delta\to0$. Using \eqref{eq: fidis}, it follows from \cite[Theorem 7.5]{billingsley2013convergence} that $Z_N$ converges weakly to $\sigma_\infty W$ on all compact intervals $[-M,M]$}. Combining this with \eqref{eq: N1/3UN} and \eqref{eq: lime} we conclude that $N^{1/3}(U_N(a)-g(a))$ is the location of the maximum of a process that weakly converges to the continuous Gaussian process
$$ \sigma_\infty(t)W(u)-\frac{|\mu'(t)|f_\infty(t)}2u^2,\ u\in\RR.$$
The above process achieves its maximum at a unique point $\mathbb{T}$ by Lemma 2.6 of \cite{KP90}, and it follows from Lemma \ref{lem: tail boundUN} that $N^{1/3}(U_N(a)-g(a))$ is uniformly tight. Hence, Corollary 5.58 in van der Vaart shows that $N^{1/3}(U_N(a)-g(a))$ converges in distribution to $\mathbb{T}$. Now, $\mathbb{T}$ is also the unique location of the maximum of the process
$$W(u)-\frac{|\mu'(t)|f_\infty(t)}{2\sigma_\infty(t)}u^2,\ u\in\RR.$$
Changing scale in the Brownian motion finally shows that
$$\left(\frac{|\mu'(t)|f_\infty(t)}{2\sigma_\infty(t)}\right)^{2/3}\mathbb{T}$$ has the same distribution as $\mathbb{Z}$, which completes the proof.
\hfill{$\Box$}

\bibliographystyle{apalike}
\bibliography{AG} 

\def\noopsort#1{}
\begin{thebibliography}{}

\bibitem[Anevski et~al., 2006]{anevski2006general}
Anevski, D., H{\"o}ssjer, O., et~al. (2006).
\newblock A general asymptotic scheme for inference under order restrictions.
\newblock {\em The Annals of Statistics}, 34(4):1874--1930.

\bibitem[Azadbakhsh et~al., 2014]{azadbakhsh2014computing}
Azadbakhsh, M., Jankowski, H., and Gao, X. (2014).
\newblock Computing confidence intervals for log-concave densities.
\newblock {\em Computational Statistics \& Data Analysis}, 75:248--264.

\bibitem[Banerjee et~al., 2018]{banerjee2018divide}
Banerjee, M., Durot, C., and Sen, B. (2018).
\newblock Divide and conquer in non-standard problems and the super-efficiency
  phenomenon.
\newblock {\em to appear in Annals of Statistics}.

\bibitem[Battey et~al., 2015]{battey2015distributed}
Battey, H., Fan, J., Liu, H., Lu, J., and Zhu, Z. (2015).
\newblock Distributed estimation and inference with statistical guarantees.
\newblock {\em arXiv preprint arXiv:1509.05457}.

\bibitem[Billingsley, 2013]{billingsley2013convergence}
Billingsley, P. (2013).
\newblock {\em Convergence of probability measures}.
\newblock John Wiley \& Sons.

\bibitem[Groeneboom et~al., 2010]{groeneboom2010maximum}
Groeneboom, P., Jongbloed, G., Witte, B.~I., et~al. (2010).
\newblock Maximum smoothed likelihood estimation and smoothed maximum
  likelihood estimation in the current status model.
\newblock {\em The Annals of Statistics}, 38(1):352--387.

\bibitem[Groeneboom and Wellner, 2001]{GW01}
Groeneboom, P. and Wellner, J.~A. (2001).
\newblock Computing {C}hernoff's distribution.
\newblock {\em J. Comput. Graph. Statist.}, 10(2):388--400.

\bibitem[Hsieh et~al., 2014]{hsieh2014divide}
Hsieh, C.-J., Si, S., and Dhillon, I. (2014).
\newblock A divide-and-conquer solver for kernel support vector machines.
\newblock In {\em International Conference on Machine Learning}, pages
  566--574.

\bibitem[Kim and Pollard, 1990]{KP90}
Kim, J. and Pollard, D. (1990).
\newblock Cube root asymptotics.
\newblock {\em Annals of Statistics}, 18:191--219.

\bibitem[Li et~al., 2013]{LiEtAl13}
Li, R., Lin, D.~K., and Li, B. (2013).
\newblock Statistical inference in massive data sets.
\newblock {\em Applied Stochastic Models in Business and Industry},
  29(5):399--409.

\bibitem[Mammen, 1991]{Mammen91}
Mammen, E. (1991).
\newblock Nonparametric regression under qualitative smoothness assumptions.
\newblock {\em Ann. Statist.}, 19(2):741--759.

\bibitem[Massart, 1990]{massart1990tight}
Massart, P. (1990).
\newblock The tight constant in the dvoretzky-kiefer-wolfowitz inequality.
\newblock {\em The Annals of Probability}, pages 1269--1283.

\bibitem[Mukerjee, 1988]{mukerjee1988monotone}
Mukerjee, H. (1988).
\newblock Monotone nonparametric regression.
\newblock {\em The Annals of Statistics}, pages 741--750.

\bibitem[Robertson et~al., 1988]{RWD88}
Robertson, T., Wright, F.~T., and Dykstra, R.~L. (1988).
\newblock {\em Order restricted statistical inference}.
\newblock Wiley Series in Probability and Mathematical Statistics: Probability
  and Mathematical Statistics. John Wiley \& Sons Ltd., Chichester.

\bibitem[Rosenthal, 1970]{rosenthal1970subspaces}
Rosenthal, H.~P. (1970).
\newblock On the subspaces of {$L_p, (p> 2)$} spanned by sequences of
  independent random variables.
\newblock {\em Israel Journal of Mathematics}, 8(3):273--303.

\bibitem[Shang and Cheng, 2017]{shang2017computational}
Shang, Z. and Cheng, G. (2017).
\newblock Computational limits of a distributed algorithm for smoothing spline.
\newblock {\em The Journal of Machine Learning Research}, 18(1):3809--3845.

\bibitem[Shi et~al., 2017]{shi2017massive}
Shi, C., Lu, W., and Song, R. (2017).
\newblock A massive data framework for m-estimators with cubic-rate.
\newblock {\em Journal of the American Statistical Association},
  (just-accepted).

\bibitem[Tang et~al., 2012]{tang2012likelihood}
Tang, R., Banerjee, M., Kosorok, M.~R., et~al. (2012).
\newblock Likelihood based inference for current status data on a grid: A
  boundary phenomenon and an adaptive inference procedure.
\newblock {\em The Annals of Statistics}, 40(1):45--72.

\bibitem[Van Der~Vaart and Van Der~Laan, 2003]{van2003smooth}
Van Der~Vaart, A. and Van Der~Laan, M. (2003).
\newblock Smooth estimation of a monotone density.
\newblock {\em Statistics: A Journal of Theoretical and Applied Statistics},
  37(3):189--203.

\bibitem[Volgushev et~al., 2017]{volgushev2017distributed}
Volgushev, S., Chao, S.-K., and Cheng, G. (2017).
\newblock Distributed inference for quantile regression processes.
\newblock {\em arXiv preprint arXiv:1701.06088}.

\bibitem[Zhang et~al., 2001]{zhang2001asymptotic}
Zhang, R., Kim, J., and Woodroofe, M. (2001).
\newblock Asymptotic analysis of isotonic estimation for grouped data.
\newblock {\em Journal of statistical planning and inference},
  98(1-2):107--117.

\bibitem[Zhang et~al., 2013]{Zhang13}
Zhang, Y., Duchi, J., and Wainwright, M. (2013).
\newblock Divide and conquer kernel ridge regression.
\newblock In {\em Conference on Learning Theory}, pages 592--617.

\bibitem[Zhao et~al., 2014]{Zhao14}
Zhao, T., Cheng, G., and Liu, H. (2014).
\newblock A partially linear framework for massive heterogeneous data.
\newblock {\em arXiv preprint arXiv:1410.8570}.

\end{thebibliography}

\end{document}